\documentclass[amscd,amssymb,verbatim,11pt]{amsart}
\numberwithin{equation}{section}
\usepackage{ulem}
\usepackage{epsfig}
\usepackage{multido}
\usepackage{color}
\usepackage{pstricks}
\usepackage{pst-all}
\usepackage{pst-math}
\usepackage{pst-func}
\usepackage{pst-3dplot}
\usepackage{enumerate}

\newtheorem{thm}{Theorem}[section]
\newtheorem{cor}[thm]{Corollary}
\newtheorem{lem}[thm]{Lemma}

\theoremstyle{definition}

\newtheorem{rem}[thm]{Remark}

\newif\ifShowLabels
\ShowLabelstrue
\newdimen\theight
\def\TeXref#1{
     \leavevmode\vadjust{\setbox0=\hbox{{\tt
            \quad\quad  {\small  \bf #1}}}%
     \theight=\ht0
     \advance\theight  by  \dp0
     \advance\theight  by  \lineskip
     \kern -\theight \vbox  to
     \theight{\rightline{\rlap{\box0}}%
      \vss}%
      }}%

    {\begin{thm}\label{#1} \ifShowLabels \TeXref{#1} \fi}%
    {\end{thm}}

    {\begin{def}\label{#1} \ifShowLabels \TeXref{#1} \fi}%
    {\end{def}}

    {\begin{lem}\label{#1} \ifShowLabels \TeXref{#1} \fi}%
    {\end{lem}}

    {\begin{cor}\label{#1} \ifShowLabels \TeXref{#1} \fi}%
    {\end{cor}}

\newcommand{\eqRef}[1]%
     {\ifShowLabels \TeXref{#1} \fi
      \begin{equation}\label{#1} }

\ShowLabelsfalse

\newcommand{\vsp}{\vskip 1em}

\newcommand{\NI}{\noindent}
\newcommand{\bea}{\begin{eqnarray}}
\newcommand{\eea}{\end{eqnarray}}

\newcommand{\be}{\begin{equation}}
\newcommand{\ee}{\end{equation}}
\newcommand{\ben}{\begin{eqnarray*}}
\newcommand{\een}{\end{eqnarray*}}
\newcommand{\lm}{\lambda}
\newcommand{\Om}{\Omega}

\newcommand{\p}{\partial}
\newcommand{\al}{\alpha}
\newcommand{\bt}{\beta}

\newcommand{\g}{\gamma}
\newcommand{\vep}{\varepsilon}
\newcommand{\dl}{\delta}
\newcommand{\D}{\Delta}

\newcommand{\lf}{\left }

\newcommand{\tx}{\text}

\newcommand{\lam}{\lambda}

\newcommand{\dt}{\delta}

\newcommand{\s}{\sigma}
\newcommand{\ep}{\varepsilon}
\newcommand{\su}{\subseteq}
\newcommand{\rar}{\rightarrow}

\newcommand{\R}{\mathbb{R}}

\newcommand{\vp}{\varphi}

\newcommand{\ve}{\varepsilon}

\newcommand{\Df}{\D_\infty}
\newcommand{\epf}{\hfill $\Box$\\}

\newcommand{\nin}{\noindent}

\newcommand{\pf}{\noindent {\bf Proof}:$\;\;$}

\title[{\large $\infty$}-Laplacian]{Inhomogeneous Dirichlet problems  involving the infinity-Laplacian}

\author[Bhattacharya, Mohammed]{Tilak Bhattacharya and Ahmed Mohammed}

\address{Department of Mathematics and Computer Science,
Western Kentucky University,
 Bowling Green, KY 42101}
\email{tilak.bhattacharya@wku.edu}

\address{Department of Mathematical Sciences, Ball State University, Muncie, IN 47306, USA}
\email{amohammed@bsu.edu}

\begin{document}

\maketitle

\begin{abstract}
Our purpose in this paper is to provide a self contained account of the inhomogeneous Dirichlet problem $\Df u=f(x,u)$ where $u$ takes a prescribed continuous data on the boundary of bounded domains. We employ a combination of Perron's method and a priori estimates to give general sufficient conditions on the right hand side $f$ that would ensure existence of viscosity solutions to the Dirichlet problem. Examples show that these sufficient conditions may not be relaxed. We also identify a class of inhomogeneous terms for which the corresponding Dirichlet problem has no solution in any domain with large in-radius. Several results, which are of independent interest, are developed to build towards the main results. The existence theorems provide substantial improvement of previous results, including our earlier results \cite{BMO} on this topic.\\

\NI{\bf Keywords.} Infinity Laplacian, Dirichlet problem, nonlinear boundary value problem, comparison principle.\\

\NI{\bf 2010 Mathematics Subject Classification.} 35J60, 35J70

\end{abstract}

\section{\bf Introduction}

\NI In this work, we continue our discussion of the Dirichlet problem for equations involving the infinity Laplacian and having inhomogeneous right hand sides. This was first initiated in \cite{BMO}, and although we stated results about existence and uniqueness of solutions applicable to various situations, it was not clear if results that allowed greater generality could be proved. Our intention in the current work is to address this matter and provide statements that apply to more general situations. To make our discussion more precise, we introduce the problem that concerns our work. Let $\Om\su\R^N,\;N\ge 2$ denote a bounded domain and $f:\Om\times \R\rar\R$ be continuous. Given $b\in C(\p\Om)$, we consider the Dirichlet problem
\eqRef{dp}
\left\{\begin{array}{rcll}
         \Df u &=& f(x,u)\;\;\;&\mbox{in}\;\Om \\[.2cm]
         u & =&b\;\;\;&\mbox{on}\;\p\Om.
       \end{array}
 \right.
\end{equation}
Here,
$$\Df u:=\langle D^2u Du,Du\rangle$$ is the infinity Laplacian. This elliptic operator is nonlinear and highly degenerate and in general does not have smooth solutions.
Let $\overline{\Om}$, denote the closure of $\Om$. By a solution to the Problem (\ref{dp}), we mean a function $u\in C(\overline{\Om})$ that satisfies (\ref{dp}) in the viscosity sense, see Section 2 for definition. In general, solutions are not known to be any better than $C^{1,\al}_{loc}$ and classical solutions, if any, form a small class. In this regard we mention the papers \cite{ESA,OS}. In the work \cite{ESA}, Evans and Savin show that infinity-harmonic functions in $\Om\su\R^2$ are $C^{1,\al}_{loc}$. As of yet, it is not known whether this is true in dimensions $N\geq 3$.  On a positive side, Evans and Smart recently provided the first breakthrough in \cite{ES} where they proved that infinity harmonic functions in $\Om\su\R^N$ are differentiable everywhere in $\Om$. For further information and motivation for studying this operator, we direct the reader to \cite{AMJ, BDM, CRM}.\\

\NI Problem (\ref{dp}) has been investigated by several authors when the right hand side $f(x,t)$ is independent of $t$. We point the reader to the recent papers \cite{AMS, LUW, LW2, LW3}. Motivated by these papers, the present authors initiated the study of the above Dirichlet problem in the paper \cite{BMO}. \\

For the reader's benefit, we start by highlighting the main results of the paper \cite{BMO}. To the best of our knowledge, \cite{BMO} is the first paper that addresses the Dirichlet Problem (\ref{dp}) in which the inhomogeneous term $f$ depends on both the variables $x$ and $t$. The paper studies Problem (\ref{dp}) when $f$ is either
non-negative or non-positive in $\Om\times \R$, under the assumption that $f(x,t)$ is bounded in $x\in\Om$ for each $t$. One of the main results in this paper shows that if $f(x,t)$ is non-decreasing in $t$, for each $x\in\Om$, then Problem (\ref{dp}) admits a solution $u\in C(\overline{\Om})$, in the viscosity sense. In contrast, it is also shown in \cite{BMO} that when $f(x,t)$ is non-increasing in $t$, for each $x\in\Om$, then solutions exist when the underlying domain $\Om$ satisfies certain size restrictions, depending on $f$ and $b$. That this is not a mere technical restriction is highlighted by an example of a non-positive decreasing inhomogeneous term $f(t)$ for which the corresponding Dirichlet Problem (\ref{dp}), for a large enough $\Om$,  fails to have a solution in $C(\overline{\Om})$.\\

Our main goal in the current work is to remove the sign and the monotonicity  restrictions, and present fairly general sufficient conditions on $f$ to ensure the existence of viscosity solutions to the Problem (\ref{dp}). We pay particular attention to developing sufficient conditions on the inhomogeneous term $f$ that would imply the existence of solution in any bounded domain. We also provide explicit conditions on $f$ that obstruct the existence of solutions to (\ref{dp}) in domains with large in-radius. Specific classes  of inhomogeneous terms are provided to exemplify the general results on non-existence. In summary, our goal, in the current work, is to explore these issues in more detail, and to present results that clarify, to some extent, the connection between the conditions that lead to existence and the those that lead to non-existence. \\

We describe in some detail the main issues that we address in this work.
The primary tool used in achieving one of the main goals of this work, namely the existence result, is a combination of an adaptation of the Perron method and a priori supremum bounds. The Perron method  approach was used in \cite{BMO}, however, its adaptation in our current work frees us of the various restrictions used in \cite{BMO} by assuming the existence of appropriate sub-solutions and super-solutions. On the face of it, this may seem to be restrictive, but as shown in this work, this actually permits greater generality. As a consequence, a great part of this work is devoted to investigating various situations that permit us to construct appropriate sub-solutions and super-solutions, see Theorem 3.1 in Section 3. As shown, the cases discussed in \cite{BMO} now become special instances of our main existence theorem of the current work. As a part of our investigations, we also look at the class of $f$ for which we can derive a priori supremum bounds. It is shown that for such $f$'s we can construct appropriate sub-solutions and super-solutions, thus leading to existence. This is proven without any restrictions on the domain $\Om$, see Theorems 5.3 and 5.5 in Section 5. We have also revisited the non-existence result in \cite{BMO} and provide sufficient conditions on $f$ that would ensure non-existence of solutions in domains of large in-radius, see Theorem 4.1 in Section 4. Thus, for arbitrary $f$, restrictions on the size of domain $\Om$ are required for existence. These assist in the construction of an appropriate sub-solution and an appropriate super-solution, and our main existence result then implies the existence of solutions to (\ref{dp}), see Theorem 3.4 in Section 3. We have also included a proof of the Lipschitz continuity of solutions to (\ref{dp}) with non-trivial right hand sides, see Theorem 2.4 in Section 2. A proof of the same result may be found in \cite{LW2}, we have provided a different proof of the same result. Also included is a Harnack inequality, which we hope will be of some independent interest, see Theorem 7.1 in Section 7.\\

Regarding the matter of uniqueness, we recall that when $f(x,t)>0$ or $f(x,t)<0$ in $\Om\times\R$, and $f(x,t)$ is non-decreasing in $t$ for each $x\in\Om$, it is known that Problem (\ref{dp}) has a unique solution. We refer the reader to \cite{BMO} for more discussion on this. As shown in \cite{LUW}, uniqueness does not hold in general if $f$ changes sign. In this work, we do not address the question of uniqueness per se, except in a special case when $f$ is a continuous function that depends on $t$ only, is non-decreasing and the boundary data is a constant. However, when $f$ is non-increasing we do not expect uniqueness to hold, in general. This is borne out by many of the examples discussed here. As a matter of fact, we have included an example, to which our main existence result applies, and multiple solutions (that are not scalar multiples) exist. It would be interesting to know what general conditions on $f(x,t)$ would lead to uniqueness of solutions to the Dirichlet problem (\ref{dp}). In the latter part of this work, we discuss some comparison principles, in this context, that apply in special situations, and may lead to uniqueness. We have also not addressed, in any great detail, circumstances that may lead to multiplicity of solutions.\\

Now a word about the layout of the paper. In Section 2, we introduce notations and definitions that will be used throughout this paper. We recall a useful comparison principle, and also prove the local Lipschitz continuity of solutions to the Dirichlet Problem (\ref{dp}) with bounded right hand sides. In Section 3, we prove our main existence result, Theorem \ref{sst}. As the first application, this theorem is used to show the existence of solutions to the Dirichlet problem (\ref{dp}) when $f(x,t)$ is bounded or is non-decreasing in $t$, for each $x\in\Om$.  Next, we apply Theorem \ref{sst} to arbitrary $f$, under size restrictions on $\Om$, to conclude existence of solutions. In Section 4, we discuss conditions under which solutions may fail to exist. In particular, we provide a sufficient condition on $f$ for which (\ref{dp}) fails to have a solution, when the underlying domain is large. We also show that if this condition fails to hold (with a slight modification), the Dirichlet problem admits a solution in any domain. Section 5 is devoted to developing a priori $L^\infty$ bounds for solutions of (\ref{dp}) when the inhomogeneous term satisfies appropriate growth conditions at infinity. These a priori bounds are then exploited to show existence of a solution. In the event these conditions are not met, a priori bounds may not hold. This is delineated by an example in Section 6. Some related issues are also discussed in the same section. In Section 7, we have included a Harnack's inequality for non-negative solutions of $\Df u\leq h(x)$ when $h$ is continuous and bounded. This generalizes the known result for non-negative infinity super-harmonic functions, and may be of some independent interest.  This section also contains a version of a comparison principle that provides estimates for solutions to (\ref{dp}), when $f(x,t)=f(t)$ is non-decreasing. In particular, it leads to uniqueness of solutions of the corresponding Dirichlet problems with constant boundary data. Finally, we have included an Appendix where we prove several lemmas required in Sections 4 and 5.
\\
\\

\section{\bf Preliminaries}

In this section, we fix some notations and state definitions that will be used throughout the paper. For easy reference, we will recall a comparison principle, that applies to the partial differential equation (henceforth referred to as PDE) in (\ref{dp}), when $f(x,t)$ is independent of $t$. We will also present a proof of the local Lipschitz continuity of solutions of (\ref{dp}). \\

We start this section by fixing some notations that will be used throughout this paper. We will work in $\mathbb{R}^N, N\ge 2$, and $\Om\subset \mathbb{R}^N$ will always stand for a bounded domain. The letters $x,\;y$ and $z$ will often denote points in $\mathbb{R}^N$, and $o$ will stand for the origin. Given $r>0$, and $x\in\mathbb{R}^N$, we will use $B_r(x)$ to indicate a ball of radius $r$ centered at $x$. For a bounded domain $\Om\su\mathbb{R}^N$, we adopt the standard notations $\overline{\Om}$ and $\p\Om$ to denote the closure and the boundary of $\Om$, respectively. We shall use $C(\Om)$ to indicate the class of functions that are continuous on $\Om$, while $C^2(\Om)$ will mean the set of twice continuously differentiable functions on $\Om$. The little oh and big oh notations will also be used. More explicitly, we will write $k(x)=O(g(x))$ as $x\rar x_0$ if there are positive constants $M$  and $\dt$ such that $|k(x)|\leq M|g(x)|$ whenever $|x-x_0|<\dt$. Likewise we use the notation $k(x)=o(g(x))$ as $x\rar x_0$ if $k(x)/g(x)\rar0$ as $x\rar x_0$. Finally, through out this work, the symbol $\s$ will denote the constant
$$\s=\frac{3^{4/3}}{4}.$$
$\;$\\
We will also have occasion to use the terms out-radius and  in-radius of a bounded domain $\Om$. By an out-ball of $\Om$, we will mean the smallest ball that contains $\Om$. The out-radius of $\Om$ is the radius of its out-ball. Similarly, an in-ball of $\Om$ is the largest ball that is contained in $\Om$, and the in-radius of $\Om$ is the radius of its in-ball.  \\

We now state the concept of viscosity solution that was introduced by Crandall, Evans and Lions \cite{CEL,CIL}.\\

\NI{\bf Definitions.} We say that a function $u\in C(\Om)$ is a sub-solution of the PDE $\Delta_{\infty}u=f(x,u)$ if for every $\psi\in C^2(\Om)$, with the property that $u-\psi$ has a local maximum at some $x_0\in \Om$, then
$$\Delta_{\infty} \psi(x_0)\ge f(x_0,u(x_0)).$$
We say that a function $u\in C(\Om)$ is a super-solution of the PDE $\Delta_{\infty}u=f(x,u)$ if for every $\psi\in C^2(\Om)$, with the property that $u-\psi$ has a local minimum at some $x_0\in \Om$, then
$$\Delta_{\infty} \psi(x_0)\le f(x_0,u(x_0)).$$
A function $u\in C(\Om)$ is a viscosity solution of $\Delta_{\infty}u=f(x,u)$ if it is both a sub-solution and a super-solution.\\

\NI In reference to the Dirichlet problem (\ref{dp}), we say $w\in C(\overline{\Om})$ is a sub-solution of Problem (\ref{dp}) if
\eqRef{ss1}
\Df w\ge f(x,w)\;\;\mbox{in $\Om$ with $w\leq b$ on the boundary $\p\Om$}.
\ee
Similarly, we say $v\in C(\overline{\Om})$ is a super-solution of Problem (\ref{dp}) if
\eqRef{ss2}
\Df v\le f(x,v)\;\;\mbox{in $\Om$ with $v\geq b$ on the boundary $\p\Om$.}
\ee$\;$\\

\begin{rem}\label{minmax} Let $\al,\bt\in C(\overline{\Om})$. The two statements below follow easily from the definition of sub (super)-solution of (\ref{dp}).

\begin{enumerate}[(i)]
\item If $\al$ and $\bt$ are sub-solutions of (\ref{dp}), then $\max\{\al,\bt\}$ is a sub-solution of (\ref{dp}).
\item If $\al$ and $\bt$ are super-solutions of (\ref{dp}), then $\min\{\al,\bt\}$ is a super-solution of (\ref{dp}).

\end{enumerate}
\end{rem}

\NI  The following versions of the comparison principle will prove useful in our subsequent discussions.
  \\
\begin{lem}[Comparison Principle]\label{luw} Suppose $f\in C(\Om)$ is and $f>0,\,f<0$ or $f\equiv0$ in $\Om$. Let $u,v\in C(\overline{\Om})$ satisfy $\Df u\geq f(x)$ in $\Om$, and $\Df v\leq f(x)$ in $\Om$.
Then
$$\sup_{\Om}(u-v)=\sup_{\p\Om}(u-v).$$
\end{lem}
\nin See \cite{AMS, BMO, J, LUW, LW3}. \\
We will often use Lemma \ref{luw} when the function $f$ is a constant function.\\

We will also need the following comparison principle which is a special case of Lemma 4.1 in \cite{BMO}.

\begin{lem}\label{mcomp}
  Suppose $f_i:\Om\rar\mathbb{R}$ are
  continuous for $i=1,2$. Let $u, v\in C(\overline{\Om})$ satisfy
  $\Delta_{\infty}u\ge f_1(x)$ and $\Delta_{\infty}v\le
  f_2(x)$. Suppose further  $f_1(x)>f_2(x)\;\;\forall\;x\in
  \Om$. Then
$$\sup_{\Om}(u-v)=\sup_{\p\Om}(u-v).$$
\end{lem}
\vsp

The next lemma establishes the local Lipschitz continuity of a sub-solution or super-solution to $\Df u=h(x)$ in $\Om$ where $h\in L^\infty(\Om)$. This result will be used throughout this work, and is an extension of a well known result, \cite{AMJ,BHA,CRM,CEG}. A similar result appears in the recent paper \cite{LW2}. Here we give a self contained proof. See also \cite{AMJ,BHA,CRM,CEG} and Lemma 3.4 in \cite{BMO}.\\

\begin{thm}\label{lpc}(Lipschitz Continuity)
Let $\al$ be a constant. Any solution $u\in C(\Om)\cap L^{\infty}(\Om)$ of

$$\Df u(x)  \geq  \al\;\;\mbox{in}\;\; \Om.$$
is locally Lipschitz continuous in $\Om$. More specifically, given $x_0\in\Om$ there is a constant $C$ that depends on $x_0,\,\mbox{diam}(\Om),\;|\al|$ and $\|u\|_{L^\infty(\Om)}$ such that
$$|u(x)-u(y)|\leq C|x-y|,\;\;\;\;\;x,y\in B_{r(x_0)/3}(x_0),$$
where $r(x_0):=\mbox{dist}(x_0,\p\Om)$. A similar result holds if $\Df u\le \al$ in $\Om$.
\end{thm}

\pf For each $x\in \Om$, let $r(x)=\mbox{dist}(x,\partial \Om).$  Set $M:=\max_{\overline{\Om}}u$, and $m:=\min_{\overline{\Om}}u$. Let us fix $x_0\in\Om$ and take
\eqRef{k}
k(x_0):=\frac{2(M-m)}{r(x_0)}+1+|\al|\mbox{diam}(\Om).
\end{equation}

For a fixed, but arbitrary $y\in B_{r(x_0)/3}(x_0)$ consider the function
$$\psi(x):=u(y)+ k|x-y|-\frac{|\al|}{2}|x-y|^2,$$ where $k:=k(x_0)$ is chosen as in (\ref{k}).

\noindent Note that $\psi \in C^\infty(\mathbb{R}^n-\{y\})$.
For $x\not=y$, we get
$$\Df\psi(x)=-|\al|(k-|\al||x-y|)^2.$$
As $k\geq 1+|\al|\mbox{diam}(\Om)$, we see that $\Df\psi\leq \al$ in $\Om\setminus\{y\}$. Taking (\ref{k}) into consideration and noting that $r(y)\geq 2r(x_0)/3$, we observe that for $x\in\p B(y,r(y))$
\begin{align*}
\psi(x)&=u(y)+kr(y)-\frac {|\al|}{2} r(y)^2\\&\geq
m+\frac{r(x_0)}{2}\left(k-\frac {|\al|}{2} r(y)\right)\\&\geq
m+\frac{r(x_0)}{2}\left(k-\frac  {|\al|}{2} \mbox{diam}(\Om)\right)\ge M\geq u(x).
\end{align*}
Therefore $u\leq \psi$ on $\p(B_{r(y)}(y)\setminus\{y\})$,
$$\Df \psi\leq \al,\;\;\mbox{and}\;\;\Df u\geq \al\;\;\;\;\mbox{in}\;\;\;B_{r(y)}(y)\setminus\{y\}.$$ Thus, the comparison principle in Lemma \ref{luw}, shows that
$u\leq \psi$ on $B_{r(y)}(y)$. Hence for any $y\in B_{r(x_0)/3}(x_0)$ and any $z\in B_{r(y)}(y)$, we have shown that
\eqRef{uu}
u(z)\leq u(y)+k|z-y|-\frac  {|\al|}{2}|z-y|^2.
\end{equation}
Now let us observe that for any $p\in B_{r(x_0)/3}(x_0)$, we have $\,B_{r(x_0)/3}(x_0)\su B_{r(p)}(p).$
Employing this observation together with (\ref{uu}), if $x,y\in B_{r(x_0)/3}(x_0)$ then we conclude that
\begin{align*} u(y)&\leq u(x)+k|x-y|-\frac {|\al|}{2}|x-y|^2,\;\;\;\;\;\;\;\;\;\;\;\mbox{and}\\
u(x)&\leq u(y)+k|x-y|-\frac  {|\al|}{2}|x-y|^2.
\end{align*} That is
$$|u(x)-u(y)|\leq \left(k-\frac  {|\al|}{2} |x-y|\right)|x-y|
\leq k|x-y|,\;\;\;\;\;\;x,y\in B_{r(x_0)/3}(x_0).$$
In conclusion, given $x_0\in\Om$, we have proved that $|u(x)-u(y)|\leq C|x-y|,\;\forall\;x,y\in B_{r(x_0)/3}(x_0)$,
where $C$ depends on $x_0, \,\mbox{diam}(\Om),\,|\al|$ and $\|u\|_{L^\infty(\Om)}$. \epf\\

\begin{cor} If $u\in C(\overline{\Om})$ is a viscosity solution (sub-solution or super-solution) of
$$\Df u(x)=h(x)\;\;\;\;\;\;\;\;\;\;(x\in\Om),$$
where $h\in C(\Om)\cap L^{\infty}(\Om)$, then $u$ is locally Lipschitz continuous in $\Om$.
\end{cor}
$\\$

\section{\bf Existence via sub-solution and super-solution method}

In this section we discuss existence of solutions to the Dirichlet problem (\ref{dp}). We place no overt restrictions on the right hand side function $f$, and the main tool used to achieve existence is the Perron method, see Theorem \ref{sst}. We will, however, assume the existence of appropriate sub-solutions and super-solutions to (\ref{dp}). The Perron method, as a tool to show the existence of viscosity solutions to fully non-linear elliptic equations, was pioneered by Hitoshi Ishii \cite{ISH1,ISH2, ISH3}. However none of these cited works seems to address the statement we seek to prove in this section, namely Theorem \ref{sst} below. As has been observed in \cite{BMO}, the existence of a solution is not guaranteed, in general, for any $f$ and any $\Om$.
However, if $\Om$ satisfies a certain size restriction, depending on $f$ and $b$, then existence will follow (see (\ref{sp})). The approach to the question of existence, in this and subsequent sections, will be to see when appropriate sub-solutions and super-solutions can be constructed so that Theorem \ref{sst} can be applied. Our effort along the way will be to try and identify the class of the functions $f$ for which this is achievable. For instance, we will discuss the case when $f$ is bounded, or when $f(x,t)$ is non-decreasing in $t$ for each $x$, and show existence of solutions of Problem (\ref{dp}).\\

\nin We now proceed to stating and proving the main result of this section. The proof is an adaptation of the proof of Theorem 3.5 in \cite{BMO}, but applies to more general situations. We assume that $f$ satisfies the following condition: for every compact interval $I\su\R$,

\eqRef{sup}
\sup_{\Om\times I}|f(x,t)|<\infty.
\end{equation}
Also recall (\ref{ss1}) and (\ref{ss2}) for the definitions of a sub-solution and a super-solution to (\ref{dp}). The following theorem holds.

\begin{thm}\label{sst} Let $f\in C(\Om\times\R,\R)$ satisfy condition (\ref{sup}). Suppose that $u_*\in C(\overline{\Om})$ is a sub-solution of (\ref{dp}) in $\Om$, and $u^*\in C(\overline{\Om})$ is a super-solution of (\ref{dp}) in $\Om$. If $u_*\leq u^*$ in $\Om$ then problem (\ref{dp}) admits a solution $u\in C(\overline{\Om})$ such that $u_*\leq u\leq u^* $ in $\Om$.
\end{thm}

\pf Let
$$\aleph_{\geq}:=\{\al\in C(\overline{\Om}): \Df \al\geq f(x,\al)\;\mbox{in}\;\Om,\;\al\leq u^*\;\mbox{in $\Om,$ and}\;\;\al\leq b\;\;\mbox{on}\;\;\p\Om\}.$$
Note that $u_*\in\aleph_{\geq}$ and $\al\leq \max_\Om u^*$ in $\Om$ for all $\al\in\aleph_{\geq}$. Therefore the set $\aleph_{\geq}$ is non-empty and bounded above.

\nin Let $\vartheta_*:=\min_\Om u_*\;\;\mbox{and}\;\;\vartheta^*:=\max_\Om u^*,$ and set
$$u(x):=\sup_{\al\in\aleph_{\geq}}\al(x),\;\;\;\;\;\;x\in\overline{\Om}.$$
It is clear from the definition that
$$\vartheta_*\leq u_*\leq u\leq u^*\leq\vartheta^* \;\;\;\tx{in} \;\;\overline{\Om}.$$
By Remark \ref{minmax} we note that $\max\{\al,u_*\}\in\aleph_{\geq}$ for any $\al\in \aleph_{\geq}$. Therefore one may also write
$$u(x):=\sup_{\al\in\widehat{\aleph}_{\geq}}\al(x),\;\;\;\;\;\;x\in\overline{\Om},$$
where
$\widehat{\aleph}_{\geq}:=\{\al\in\aleph_{\geq}:\al\geq u_*\;\;\tx{in}\;\Om\}$.\\

\nin First we show that $u\in C(\overline{\Om})$ and $u=b$ on $\p\Om$. Let us pick $C_*\leq 0$ and $C^*\geq 0$ such that $$C_*\;\leq \inf_{\Om\times[\vartheta_*\,,\,\vartheta^*]}f(x,t)\;\leq
\;\sup_{\Om\times[\vartheta_*\,,\,\vartheta^*]}f(x,t)\leq \;C^*.$$
Then $\Df\al\geq C_*$ for all $\al\in\widehat{\aleph}_{\geq}$. Therefore by Theorem \ref{lpc}, $u$ is locally Lipschitz in $\Om$. We now proceed to show that $u$ is continuous on $\partial \Om$, and $u=b$ on $\p\Om$. Let $\al_*,\;\bt^*\in C(\overline{\Om})$ such that
\begin{align*}
\Df \al_*&=C^*\;\;\mbox{in}\;\Om,\;\;\;\mbox{and}\;\;
\al_*=b\;\;\;\mbox{on}\;\;\p\Om\\
\Df\bt^*&=C_*\;\;\mbox{in}\;\Om,\;\;\;
\mbox{and}\;\;\bt^*=b\;\;\;\mbox{on}\;\;\p\Om.
\end{align*}
For a proof of the existence of $\al_*$ and $\bt^*$, see \cite{BMO, LUW}.
Note that since $\Df u^*\leq C^*$ and $\al_*\leq u^*$ on $\p\Om$, by Lemma \ref{luw}, we actually have $\al_*\leq u^*$ in $\Om$,  Likewise, we have $u_*\leq \bt^*$ in $\Om$.
\nin We define
$$\al:=\max\{\al_*,u_*\}\;\;\;\mbox{and}\;\;\;\;\bt:=\min\{\bt^*,u^*\}.$$
We notice that $\al$ and $\bt$ are both in $C(\overline{\Om})$, and
$$\vartheta_*\leq u_*\leq \al,\;\bt\leq u^*\leq \vartheta^*\;\;\; \mbox{in}\;\;\Om,\;\;\tx{and}\;\;\al=\bt=b\;\;\tx{on}\;\p\Om.$$

Our goal is to show that $\al$ is a sub-solution and lies in $\widehat{\aleph}_{\geq}.$
Let $\vp\in C^2(\Om)$ such that $\al-\vp$ has a local maximum at $x_0\in\Om$. Suppose $\al(x_0)=\al_*(x_0)$. Then the following holds in a neighborhood $N$ of $x_0$.
$$\al_*-\vp\leq \al-\vp\leq \al(x_0)-\vp(x_0)=\al_*(x_0)-\vp(x_0).$$
Thus $\al_*-\vp$ has a local maximum at $x_0$, and since $\vartheta_*\leq \al(x_0)\leq \vartheta^*$, we have
$$\Df \vp(x_0)\geq C^*\geq f(x_0,\al(x_0)).$$
On the other hand, if $\al(x_0)=u_*(x_0)$, then it is easily seen that $u_*-\vp$ has a maximum at $x_0$, and hence
$$\Df \vp(x_0)\geq f(x_0,u(x_0))=f(x_0,\al(x_0)).$$
Thus in any case, we see that
$$\Df \al\geq f(x,\al)\;\;\mbox{in}\;\Om\;\;\;\mbox{and}\;\;\al=b\;\;\mbox{on}\;\;\p\Om,$$
and therefore $\al\in\hat{\aleph}_{\geq}.$\\

\nin Now, let us recall that for any $\g\in \aleph_{\geq}$, we have $\Df \g \geq C_*$ in $\Om$, and $\g\le b$ on $\p\Om$. Therefore, by Lemma \ref{luw}, we see that $\g\leq \bt^*$ in $\Om$, and hence, $u\leq \bt^*$ in $\Om$. Therefore we see that $\al\leq u\leq \bt$ in $\overline{\Om}$. Since $\al,\;\bt\in C(\overline{\Om})$, this proves the continuity of $u$ on $\p\Om$ and that $u=b$ on $\p\Om$.\\

\nin Next we show that $u$ is a viscosity sub-solution, and note this would imply that $u\in \aleph_{\geq}$. In the last step, we will show that it is also a super-solution providing us with a solution to (\ref{dp}).\\

\nin Let $\psi\in C^2(\Om)$ and $x_0\in \Om$ be such that $u-\psi$ has a local maximum at $x_0$, in other words $u(x)-\psi(x)\le u(x_0)-\psi(x_0)$ for $x$ in some small ball $B_{\rho}(x_0).$ Our goal is to show that $\Delta_{\infty}\psi(x_0)\ge f(x_0,u(x_0))$. In what follows, we fix $B_{\rho}(x_0)\subset \Om$, and $0<\dt<\rho^4$. Since $u(x_0)=\sup_{\al\in\aleph_{\geq}}\al(x_0)$, we pick a sequence $\{\al_k\}$ in $\aleph_{\geq}$ such that $u(x_0)-\al_k(x_0)<\dt/k$ for each positive integer $k$.
Clearly, for $x$ in $B_{\rho}(x_0)$,
\eqRef{tt}
\al_k(x)-\psi(x)\le u(x)-\psi(x)\le u(x_0)-\psi(x_0)\leq \al_k(x_0)-\psi(x_0)+\dl/k.
\end{equation}
Rewriting, we get
$$
\al_k(x)-\psi(x)-\dl/k\le \al_k(x_0)-\psi(x_0),\;\;x\in B_{\rho}(x_0).
$$
Thus for all $x\in B_{\rho}(x_0)\setminus \overline{B}_{\left(\dl/k\right)^{1/4}}(x_0)$ we have
$$
\al_k(x)-[\psi(x)+|x-x_0|^4]<\al_k(x)-\psi(x)-\dl/k\le \al_k(x_0)-\psi(x_0).
$$

\nin This inequality shows that the maximum of the function $\al_k(x)-[\psi(x)+|x-x_0|^4]$ on $\overline{B}_{\rho}(x_0)$, occurs at some $x_k$ in the set  $\overline{B}_{\left(\dt/k\right)^{1/4}}(x_0)$. In particular,
\eqRef{i}
\al_k(x_k)-[\psi(x_k)+|x_k-x_0|^4]\geq \al_k(x_0)-\psi(x_0).
\end{equation}
Since $\al_k\in\aleph_{\geq}$, the definition of sub-solution gives
\eqRef{dto}
\Delta_{\infty}\psi(x_k)+O(\sqrt{\dl/k})\ge f(x_k,\al(x_k)).
\end{equation}
Combining (\ref{tt}) and (\ref{i}) we see that
$$\al_k(x_0)-\psi(x_0)\leq \al_k(x_k)-[\psi(x_k)+|x_k-x_0|^4]\leq u(x_0)-\psi(x_0)-|x_k-x_0|^4.$$
The above inequalities show that $\lim_{k\rar\infty}\al_k(x_k)=u(x_0)$. Thus letting $k\rar\infty$ and using this limit in (\ref{dto}) we conclude that $$\Delta_{\infty}\psi(x_0)\ge f(x_0,u(x_0)).$$ Thus $u$ is a sub-solution. Recalling that $u=b$ on $\partial \Om$, we see $u\in\aleph_{\geq}.$\\

\noindent The final step is to show that $u$ is a super-solution. Let $\psi\in C^2(\Om)$ and $x_0\in\Om$ be such that
$u-\psi$ has a local minimum at $x_0$. Suppose that $u(x_0)=u^*(x_0)$. Then for $x$ near $x_0$, we have the following. $$u^*(x)-\psi(x)\ge u(x)-\psi(x)\geq u^*(x_0)-\psi(x_0).$$
Clearly
$\Df\psi(x_0)\leq f(x_0,u^*(x_0))=f(x_0,u(x_0))$. So suppose $u(x_0)<u^*(x_0)$, and assume that  $\Delta_{\infty}\psi(x_0)>f(x_0,u(x_0)).$  We will show that this violates the definition of $u$. We achieve this in four steps.\\
\noindent (i) Let $d(x_0)$ denote the distance of $x_0$ from $\p\Om$. Since $u-\psi$ has a local minimum at $x_0$, there is a small ball $B_{\rho}(x_0)$ such that
$$u(x)-\psi(x)\ge u(x_0)-\psi(x_0),\;\forall\;x\in B_{\rho}(x_0).$$ Define $\phi(x):=\psi(x)+(u(x_0)-\psi(x_0))$. It is easy to see that, on $B_\rho(x_0)$,
\eqRef{strc}
\phi(x_0)=u(x_0),\;\;u(x)-\phi(x)\ge 0,\;\;\mbox{and}\;\;
\Delta_{\infty}\phi(x_0)=\Delta_{\infty}\psi(x_0)>f(x_0,\phi(x_0)).
\end{equation}
On account of the continuity of $f(x,t)$ we may choose $0<\ve_0<\min\{1,\rho, (d(x_0)/2)^8\}$, small enough, such that  $$\Delta_{\infty}\phi(x_0)>f(x_0,\phi(x_0)+\ve),\;\forall\;0<\ve\le \ve_0.$$

\noindent (ii) Let $\ve_0$ be as in Step (i). For $0<\ve\le \ve_0$,  define $\phi_{\ve}(x):=\phi(x)-\sqrt{\ve}|x-x_0|^4+\ve.$
Computation shows that
$$\Delta_{\infty}\phi_{\ve}(x)=\Delta_{\infty}\phi(x)+O(\sqrt{\ve}|x-x_0|^2),
\;\;\mbox{as $x\rightarrow x_0$}.$$
From (\ref{strc}),
$$\Df \phi_\ep(x_0)=\Df \phi(x_0)>f(x_0,\phi(x_0)+\ep)= f(x_0,\phi_{\ve}(x_0)).$$
We show that there is an $\ve_1$, with $0<\ve_1\le \ve_0\,$ small, such that $\Delta_{\infty}\phi_{\ve_1}(x)>f(x,\phi_{\ve_1}(x)),\;\forall\;x\in B_{\ve_1^{1/8}}(x_0).$ To prove this, we exploit the continuity of $f$. Assume to the contrary. Then, for each $\ve>0$ with $\ve\rightarrow 0$, there is an $x_{\ve}\in B_{\ve^{1/8}}(x_0)$ such that
$\Delta_{\infty}\phi_{\ve}(x_{\ve})\le f(x_{\ve},\phi_{\ve}(x_{\ve}))$. Since $x_{\ve}\rightarrow x_0$, we observe that
$$\lim_{\ve\rightarrow 0}\Delta_{\infty}\phi_{\ve}(x_{\ve})=\Delta_{\infty}\phi(x_0)\;\;\mbox{and}\;\;\lim_{\ve\rightarrow 0}f(x_{\ve},\phi_{\ve}(x_{\ve}))=f(x_0,\phi(x_0)).$$
We conclude that $\Delta_{\infty}\phi(x_0)\le f(x_0,\phi(x_0))$. We have a contradiction and the claim is proved. \\

\noindent Since $\phi(x_0)<u^*(x_0)$ we can assume that $\ve_1$ is small enough that $\phi(x)+\ve_1\leq u^*(x)$ for all $x\in B_{\ve_1^{1/8}}(x_0)$.  Moreover, there is an $0<s_1<\ep_1^{1/8}$ such that $u(x)<\phi_{\ep_1}(x)$ for all $x\in B_{s_1}(x_0)$. Thus
\begin{equation}\label{pe}
\left\{\begin{array}{ll}
         \Delta_{\infty}\phi_{\ve_1}(x)>f(x,\phi_{\ve_1}(x)),\;\;&\forall\;x\in B_{\ve_1^{1/8}}(x_0) \\[.3cm]
         u(x)<\phi_{\ve_1}(x),
\;\;&\forall\;x\in B_{s_1}(x_0).
       \end{array}\right.
\end{equation}

\noindent (iii) We recall and summarize the conclusions obtained so far. Firstly, from Step (i)
\ben
&&u(x)\ge \phi(x),\;\;\forall\;x\in B_{\rho}(x_0),\;\mbox{and}\\
&&u(x)-\phi_{\ve_1}(x)=u(x)-\phi(x)+\sqrt{\ve_1}|x-x_0|^4-\ve_1> 0,\;\forall\;x\in B_{\rho}(x_0)\setminus \overline{B}_{\ve_1^{1/8}}(x_0).
\een
Combining the above with (\ref{pe}) we conclude:

\begin{equation}\label{ab1}
\left\{\begin{array}{lll}
        (a)\;\;&\;\; \Delta_{\infty}\phi_{\ve_1}(x)>f(x,\phi_{\ve_1}(x)),\;&
         \forall\;x\in B_{\ve_1^{1/8}}(x_0)  \\[.3cm]
         (b)\;\;&\;\;\phi_{\ve_1}(x)<u^*(x),\;&\forall\;x\in B_{\ve_1^{1/8}}(x_0)  \\[.3cm]
         (c)\;\;&\;\; u(x)<\phi_{\ve_1}(x),\;&\forall\;x\in B_{s_1}(x_0)  \\[.3cm]
         (d)\;\;&\;\; u(x)>\phi_{\ve_1}(x),\;&\forall x\in B_{\rho}(x_0)\setminus \overline{B}_{\ve_1^{1/8}}(x_0).
       \end{array}\right.
\end{equation}$\;$\\
\nin The conclusions in (\ref{ab1}) will be instrumental in obtaining a contradiction.\\

\noindent (iv) We now define
\ben
\widehat{u}(x)=\left\{\begin{array}{cl}u(x)\;\;\;&\mbox{if}\;\;x\in \Om\setminus \overline{B}_{\ve_1^{1/8}}(x_0),\\[.2cm]
\sup\{\phi_{\ve_1}(x),\;u(x)\}\;\;&\mbox{if}\;\;x\in  B_{\ve_1^{1/8}}(x_0).\end{array}\right.
\een
It is clear that $\widehat{u}\in C(\overline{\Om})$, and $u_*\leq u\leq \widehat{u}\leq u^*$ in $\Om$ with $\widehat{u}=b$ on $\p\Om$. We want to show that $\widehat{u}\in \aleph_{\geq}$. To this end we take  $\widehat{\psi}\in C^2(\Om)$ such that $\widehat{u}-\widehat{\psi}$ has a local maximum at $y\in\Om$.
To be specific, we have $\widehat{u}(x)-\widehat{\psi}(x)\leq \widehat{u}(y)-\widehat{\psi}(y)$ for $x$ in some ball $B_\dt(y)$.
Then either $\widehat{u}(y)=u(y)$ or $\widehat{u}(y)=\phi_{\ve_1}(y)$. We look at each case separately. Let us take the case
$\widehat{u}(y)=u(y)$ first. Noting that $u\leq \widehat{u}$ in $\Om$ we have, for each $x\in B_\dt(y)$,
$$u(x)-\widehat{\psi}(x)\leq\widehat{u}(x)-\widehat{\psi}(x)\leq\widehat{u}(y)-\widehat{\psi}(y)=u(y)-\widehat{\psi}(y).$$
Thus $u-\widehat{\psi}$ has a local maximum at $y$. Since $u$ is a sub-solution, it follows that $\Delta_{\infty}\widehat{\psi}(y)\ge f(y,u(y))=f(y,\widehat{u}(y))$.
Suppose now that $\widehat{u}(y)=\phi_{\ve_1}(y)$.
Without loss of generality, we can assume $u(y)<\phi_{\ve_1}(y)$. In view of (\ref{ab1})(d) this inequality implies that $y\in B_{\ve_1^{1/8}}(x_0)$.  Then, again noting that $\phi_{\ve_1}\leq \widehat{u}$, we find
$$\phi_{\ve_1}(x)-\widehat{\psi}(x)\leq \widehat{u}(x)-\widehat{\psi}(x)\leq \widehat{u}(y)-\widehat{\psi}(y)=\phi_{\ve_1}(y)-\widehat{\psi}(y),\;\;\;\;\;\;x\in B_{\ve_1^{1/8}}(x_0)\cap B_{\dl}(y),$$
and therefore
$\phi_{\ve_1}-\widehat{\psi}$ has a local maximum at $y$. This implies that $\Df \phi_{\ve_1}(y)\leq \Df\widehat{\psi}(y)$. This together with (\ref{ab1})(a) shows that $\Delta_{\infty}\widehat{\psi}(y)\ge f(y,\phi_{\ve_1}(y))=f(y,\widehat{u}(y)).$
Thus, in either case we see that $\widehat{u}$ is a sub-solution and thus lies in $\aleph_{\geq}$. However by (\ref{ab1})(a) we see that $\widehat{u}>u$ in $B_{s_1}(x_0)$, and this contradicts the definition of $u.$ Thus $u$ is a super-solution, and this completes the proof that $u$ is a solution to the Dirichlet problem (\ref{dp}) on $\Om$. $\Box$\\
\\
\begin{rem}
\NI Let $f, \;u_*$ and $u^*$ be as in the above Theorem \ref{sst}. Then $-u^*$ is a sub-solution and $-u_*$ is a super-solution of the Dirichlet problem (\ref{dp}) with $f(x,t)$ replaced by $-f(x,-t)$ and the boundary data $b$ replaced by $-b$. Let us consider another class $\aleph_{\geq}$ defined as follows.
$$\;\;\;\;\;\;\;\;\aleph_{\geq}:=\{\al\in C(\overline{\Om}): \Df \al\geq -f(x,-\al)\;\mbox{in}\;\Om,\;\;\al\leq -u_*\;\;\mbox{in $\Om\;$ and}\;\al\leq -b\;\mbox{on}\;\p\Om\}.$$
Next, define
$$\aleph_{\leq}:=\{\bt\in C(\overline{\Om}):\; \Df \bt\leq f(x,\bt)\;\mbox{in}\;\Om,\;\;\bt\geq u_*\;\;\mbox{in $\Om$, and}\;\bt\geq b\;\mbox{on}\;\p\Om\}.$$
If by $-\aleph_{\leq}$ we denote the set $\{-\bt:\;\bt\in \aleph_{\leq}\}$, then $\aleph_{\geq}=-\aleph_{\leq}$.
Therefore $w(x)=\sup_{\al\in\aleph_{\geq}}\al(x)$ is a solution of the problem (\ref{dp}) with right-hand side $-f(x,-t)$ and boundary data $-b$.
That is $v=\inf_{\bt\in\aleph_{\leq}}\bt(x)=-w$ is a solution of the Dirichlet problem (\ref{dp}). Note that the solution $u$ constructed in the above theorem belongs to $\aleph_{\leq}$, and hence $v\leq u$ in $\Om$. The function $v$ is called the minimal solution and $u$ is called the maximal solution relative to the pair $(u_*,u^*)$ of sub-super solutions to the Dirichlet problem (\ref{dp}).\\
\end{rem}

\noindent Our goal in the rest of the section is to discuss various situations to which Theorem \ref{sst} may be applied. The following corollary addresses the case when $f$ is bounded or when $f$ is non-decreasing in the second variable. A version was proven earlier in \cite{BMO}. We derive this as a consequence of Theorem \ref{sst}.
Also, as will be clear from the proof, no domain restrictions are needed in this case. \\

\begin{cor}\label{bd} Let $f\in C(\Om\times\R,\R)$ and $b\in\p\Om$. If $f(x,t)$ is non-decreasing in $t$ for each $x\in\Om$ and $\sup_\Om |f(x,t)|<\infty$ for each $t\in\R$, or $f$ is bounded in $\Om\times\R$, then problem (\ref{dp}) admits a solution $u\in C(\overline{\Om})$.
\end{cor}

\pf As the proof for the bounded case follows along similar lines, we only consider the case when $f(x,t)$ is non-decreasing in $t$ for each $x\in\Om$.  According to Theorem \ref{sst} we need only find a sub-solution $u_*$ and a super-solution $u^*$ of (\ref{dp}) in $C(\overline{\Om}).$ To construct a sub-solution, we fix a positive constant $C$ such that $C\geq \left(\sup_\Om f(x,\ell)\right)^{1/3}$, where $\ell:=\inf_{\p\Om}b$. Fix $z\in\p\Om$, and let $d$ be a constant such that $d\leq \ell/C-\s\; \mbox{diam}(\Om)^{4/3}$, and we define
$$u_*(x)=C(\s |x-z|^{4/3}+d).$$ Then $u_*\leq\ell$ in $\Om$, and $\Df u_*=C^3\geq f(x,\ell)\geq f(x,u_*)$ in $\Om$. For a super-solution, we take a positive constant $C$ such that $C\geq -\left(\inf_{\Om}f(x,L)\right)^{1/3}$, where $L=\sup_{\p\Om}b$. We now define
$$u^*(x):=C(d-\s|x-z|^{4/3})$$ where $d$ is chosen such that $d\leq L/C-\s\;\mbox{diam}(\Om)^{4/3}$. Then $u^*\geq L$ in $\Om$, and clearly $\Df u^*=-C^3\leq f(x,L)\leq f(x,u^*)$. By Lemma \ref{luw} we also have $u_*\leq u^*$ on $\Om$. We now invoke Theorem \ref{sst} to conclude the proof.\epf\\
\\
We now state a general existence result when $f(x,t)$ is just continuous on $\Om\times\R$. As shown in \cite{BMO}, one needs a restriction, in general, on the size of the underlying domain $\Om$. We will
revisit this issue in Section 4.\\

Let $\Om\su\R^N$ be a bounded domain, and $b\in C(\Om)$. Set $\ell:=\inf_{\p\Om}b$ and $L:=\sup_{\p\Om} b$. For $\eta\geq 0$, we define
\eqRef{cons}
C(\eta):= \max\left\{\left(\sup_{\Om\times[\ell-\eta,\;\ell]}
f^+(x,t)\right)^{1/3},\,-\left(\inf_{\Om\times[L,\;L+\eta]}f^-(x, t)\right)^{1/3}\right\}.
\end{equation}
The following condition relates $b$ to the domain $\Om$ and the inhomogeneous term $f$
\eqRef{sp}
\mbox{diam}(\Om)<\sup_{\eta>0}\left(\frac{\eta}{\s C(\eta)}\right)^{3/4}.
\end{equation}
The right hand side expression in (\ref{sp}) is understood to be infinity if $C(\eta)=0$ for some $\eta\geq 0$.

\begin{thm}\label{decf}
Suppose that $f\in C(\Om\times\R,\,\R)$ satisfies condition (\ref{sup}). If $\Om$ is a bounded domain for which condition (\ref{sp}) holds, then Problem (\ref{dp}) admits a solution $u\in C(\overline{\Om}).$
\end{thm}

\pf In view of Theorem \ref{sst} it is enough to construct a sub-solution $u_*\in C(\overline{\Om})$ of (\ref{dp}) and a super-solution $u^*\in C(\overline{\Om})$ of (\ref{dp}).
\nin Suppose first that $f^+(x,\ell)\equiv 0$ in $\Om$, and $f^-(x,L)\equiv 0$ in $\Om$. Then $u_*(x)\equiv\ell$ and $u^*(x)\equiv L$ would serve as the desired sub-solution and super-solution, respectively. Therefore, in this case Theorem \ref{sst} implies that Problem (\ref{dp}) admits a solution in $\Om$ for any bounded domain $\Om$.\\

\nin So we assume that $f^+(x,\ell)\not\equiv 0$ in $\Om$ or $f^-(x,L)\not\equiv 0$ in $\Om$. Then, the continuity of $f$ implies that $C(\eta)>0$ for some $\eta>0$.  We note that condition (\ref{sp}) implies
$$\mbox{diam}(\Om)\leq\left(\frac{\dt}{\s C(\dt)}\right)^{3/4}$$
for some $\dt>0$. We proceed with our construction of a sub-solution $u_*$ of Problem (\ref{dp}). For notational convenience we will use $C$ instead of the positive constant $C(\dt)$. We select a number $d$ such that
$$\frac{\ell-\dt}{C}\le d\le \frac{\ell}{C}-\s \mbox{diam}(\Om)^{4/3}.$$
Let $z\in\p\Om$ and take $u_*(x)=C(\s |x-z|^{4/3}+d).$ It is easily seen that
$\ell-\dt\le u_*\le \ell$ in
$\overline{\Om}$ and hence $\Df u_*=C^3\ge f(x,u_*)$.\\

\NI We now construct $u^*$. We fix $z\in\p\Om$. We choose a constant $d$ such that
$$\frac{L}{C}+\s\mbox{diam}(\Om)^{4/3}\leq d\leq\frac{L+\dt}{C}.$$
We define $u^*(x):=C(d-\s|x-z|^{4/3})$.  Moreover, by (\ref{cons}) and (\ref{sp}) we see that $L\leq v\leq L+\dt$ in $\overline{\Om}$, and therefore $\Df u^*=-C^3\leq f(x,u^*)$, in $\Om$. The comparison principle, Lemma \ref{luw}, implies that $u_*\leq u^*$ in $\Om$, and hence an appeal to Theorem \ref{sst} concludes the proof of the theorem.\epf
\\

\nin{\bf Example:} We present two examples to illustrate Theorem \ref{decf}.

\begin{enumerate}[(i)]

\item Let $c<d$ be real numbers, and $\al$ and $\bt$ be positive odd integers. Take $f(x,t)=a(x)(t-c)^\al(t-d)^\bt$
where $a(x)\in C(\Om)\cap L^{\infty}(\Om)$ is non-negative. Suppose that $b\in C(\p\Om)$ such that $c<\ell<d<L$, where $\ell=\inf_\Om b$ and $L=\sup_\Om b$. Then Problem (\ref{dp}) has a solution in $C(\overline{\Om})$. The same conclusion holds if $a(x)\leq 0$ in $\Om$ and $\ell<c<L<d$.\\

\item As another example let us take $f(x,t)=-e^t$. Then $C(\eta)=e^{(L+\eta)/3}.$
Thus
$$\sup_{\eta>0}\left(\frac{\eta}{\s C(\eta)}\right)^{3/4}=\frac{1}{\s^{3/4} e^{L/4}}\sup_{\eta>0} \eta^{3/4} e^{-\eta/4}=\frac{3^{3/4}}{\s^{3/4} e^{(L+3)/4}}.$$$\;$\\

\nin Therefore, if $\Om$ satisfies
$$\tx{diam}(\Om)<\left(\frac3\s\right)^{3/4}e^{-(3+L)/4},$$ then Problem (\ref{dp}) with $f(x,t)=-e^t$ and $b\in C(\p\Om)$ has a solution $u\in C(\overline{\Om})$. For instance if we take $b\equiv0$, and $\Om$ is a bounded domain such that $$\tx{diam}(\Om)<\left(\frac{3}{\s e} \right)^{3/4},$$ then the corresponding Dirichlet pronblem (\ref{dp}) admits a solution.$\;$\\

We recall that in \cite{BMO}(see Part II of Appendix) it was shown that if the in-radius $R$ of $\Om$ satisfies
$$R>\left(\frac{3}{\s}\right)^{3/4},$$ then the Dirichlet problem (\ref{dp}) has no solution $u\in C(\overline{\Om})$ when $f(x,t)=-e^t$ and $b\equiv 0$ on $\p\Om$.
\end{enumerate}
 In Section 4 we will revisit Example (ii) above to show that solutions to (\ref{dp}) may fail to exist for a large class of inhomogeneous terms $f(x,t)$ if $\Om$ has a large in-radius.   See Remark \ref{eig2}.\\

\vspace*{1cm}

\section{\bf Sufficient Conditions for Existence and Non-Existence}

In this section we discuss sufficient conditions on the right hand side term $f(x,t)$ in the Dirichlet problem (\ref{dp}) that would provide somewhat more definitive information on the existence or non-existence of a solution in $C(\overline{\Om})$. On one hand, we state a specific criteria on $f(x,t)$ that would guarantee the existence of a solution in any bounded domain $\Om\su\R^N$.  On the other hand, we provide a complementary criterion on the non-homogeneous term in (\ref{dp}) such that the problem fails to admit any solution in $C(\overline{\Om})$ when $\Om$ has a sufficiently large in-radius. This latter result places the example in \cite[Part I of Appendix]{BMO} in a more general framework. We will provide specific examples that exemplify the criterion for non-existence of a solution. We point out that our discussion is limited to solutions of (\ref{dp}) with non-positive inhomogeneous terms.\\

\nin Let $\Om$ be a bounded domain; for ease of presentation, we will use the following altered version of (\ref{dp}).
Consider the following Dirichlet problem.
\eqRef{decrr}
\Df u=-f(x,u)\;\;\mbox{in $\Om$ with $u=b$ on $\p\Om$},
\ee
where $b\in C(\p\Om, \R)$, and $f\in C(\Om\times\R,[0,\infty))$, satisfies the condition (\ref{sup}), which we recall, may be stated as
$$\sup_{\Om\times I}f(x,t)<\infty,\;\;\mbox{for any compact interval $I$}.$$
For the rest of this section, the functions $f$ and $b$ will be as described above.
\nin We observe that if we replace $u$ in (\ref{decrr}) by $v=-u$ the above becomes $\Df v=\widetilde{f}(x,v)$, where $\widetilde{f}(x,t)=f(x,-t)$ is positive. Thus our conclusions apply to both the situations, that is, when the inhomogeneous term is either non-positive or non-negative.\\

\NI Let $\ell:=\inf_{\p\Om}b$, and define
$$h(t):=\inf_{\Om \times[t,\infty)}f(x,s),\;\;\;\;\;\mbox{for every $t\ge \ell$}.$$
Clearly the function $h(t)$ is non-decreasing and $h(t)\ge 0$. We remark that the function

$$\eta(t)=\frac{1}{t-\ell}\int^t_\ell h(s)\;ds,\;\;\;\;\;t>\ell$$
is non-decreasing, continuous on $[\ell,\infty)$, and satisfies $0\le h(\ell)\leq \eta(t)\le h(t)$ for all $t\geq \ell$. Therefore, without loss of generality, we can assume in our subsequent discussion that $h$ is continuous in $[\ell,\infty)$.
We make the following two assumptions.

\begin{align}
&h(t)>0 \;\;\tx{for}\;\;t>\ell.\label{dd1}\\[.3cm]
&\sup_{a>\ell}\zeta(a)=M_f<\infty,\;\;\;\tx{where} \label{dd2}\\[.3cm]\notag\\[-.7cm]
&\;\;\;\;\;\;\;\;\;\;\;\;\;\;\;\;\zeta(a)=
\int_\ell^a\frac{1}{\sqrt[4]{H(a)-H(t)}}\,dt,\;\;\;a>\ell\;\;\;\;\tx{and}\;\;
\;\;\;H(t)=\int_\ell^t h(s)\,ds,\;\;\;\;\;t>\ell.\notag
\end{align}

We now state a non-existence result.

\begin{thm}\label{uss} Let $\Om$ be a bounded domain and $B_R(z)$, for some $z\in\Om$, be its in-ball. Assume that, $h$, as defined above, satisfies (\ref{dd1}) and (\ref{dd2}).
Then the Dirichlet problem (\ref{decrr}) has no solutions $u\in C(\overline{\Om})$, besides, possibly the constant solution, if $R>M_f/\sqrt{2}$.
\end{thm}

\pf We argue by contradiction. Let $u\in C(\overline{\Om})$ be a solution of (\ref{decrr}) in $\Om$. Then $u$ satisfies the strong minimum principle as it is infinity super-harmonic in $\Om$. Thus $u(x)>\ell$ for $x\in \Om$, moreover, Theorem \ref{lpc} shows that $u$ is locally Lipschitz continuous in $\Om$. The idea is to employ an auxiliary function that will help us in deriving an estimate involving the solution $u$ and the in-radius $R$ of $\Om$. It will then follow from this estimate that if $R$ is large enough then there are no solutions to (\ref{decrr}). \\

\NI Recall that $B_R(z)\su\Om$. For $0\le r\le R,$ define
$$m(r):=\inf_{B_r(z)}u.$$ Thus $m(r)>\ell,$ for $0\le r<R,$ and that $m(r)$ is concave in $r$. We also observe that
$f(x,u(x))\ge \inf_{\Om\times[u(x),\infty)}f(x,s)=h(u(x)).$ Furthermore, \\
\begin{enumerate}[(i)]
\item $m(r)$ is decreasing in $r$,
\item $u(x)\ge m(|x-z|)>m(r),\;x\in B_r(z),$ and
\item $f(x,u(x))\ge h(m(|x-z|))\ge h(m(r))>0,\;x\in B_r(z).$\\
\end{enumerate}
\NI Let $w\in C([0,R])$ be defined by the equation
$$w(r)=w(0)-3^{1/3}\int_0^r\left[\int_0^t h(m(s)) ds\right]^{1/3} dt,\;\;\;0\le r<R,$$
where $w(0)$ is so chosen that $w(R)=\ell$. Clearly, $w\in C^2((0,R))\cap C^1([0,R))$ and
\eqRef{h2h2}
(w^{\prime}(r))^2w^{\prime\prime}(r)=-h(m(r)),\;0<r<R,\;\;\;\mbox{and}\;\;w^{\prime}(0)=0.
\ee
Set $v(x)=w(|x-z|),\;x\in B_R(z),$ it is easily checked that
$$\Df v(x)=-h(m(|x-z|)),\;\;\;\mbox{$x\in B_R(z)\setminus\{z\}$}.$$
We now show that $v$ is a viscosity sub-solution in $B_R(z)$. Let $\psi\in C^2(B_R(z))$ be such that $v-\psi$ has a maximum at some point $y\in B_R(z)$. Since $v\in C^2$ except at $z$, it is sufficient to consider the case when $y=z$. Since $v$ is $C^1$ at $z$ we see that $0=D v(z)=D\psi(z)$ and hence
$\Df \psi(z)=0\ge -h(m(0))$ as desired. Next we show that $v$ is a super-solution. Let $\psi\in C^2(B_R(z))$ and $y\in B_R(z))$ be a point of minimum of $v-\psi$. Once again we will consider the case when $y=z$. Then
$$w(|x-z|)-w(0)\ge \langle D\psi(z),(x-z)\rangle+\frac{1}{2}\langle D^2\psi(z)(x-z),(x-z)\rangle+o(|x-z|^2),\;\;\;x\rightarrow z.$$
Since $w^{\prime}(0)=0$, we again observe that $D\psi(z)=0$. Next by taking $x=z\pm re,\;t>0$ and $e$, any unit vector, we see that
$$\langle D^2\psi(z)e,e\rangle \le \liminf_{r\rightarrow 0}\;\frac{-2\sqrt[3]{3}}{r^2}\int_0^r\left(\int_0^t h(m(s))\;ds\right)^{1/3}\;dt=-\infty,$$
since $h(m(0))>0$. This last conclusion contradicts the choice of $\psi$ as a $C^2$ function. Thus we have shown that $v$ solves, in the sense of viscosity,
\eqRef{vm}
\Df v(x)=-h(m(|x-z|))<0,\;\;\mbox{$x\in B_R(z)$ with $v=\ell$ on $|x-z|=R$}.
\ee
\\
\NI We now derive an estimate relating $u,\;v$ and $R$. Firstly, by the observations in (i), (ii) and (iii), we note that
$$\Df u(x)=-f(x,u(x))\le -h(m(|x-z|))<0,\;\mbox{for $x\in B_R(z),$ with $u\ge \ell$ on $\p B_R(z).$}$$
Recalling (\ref{vm}) and the comparison principle, Lemma \ref{luw}, it follows that $u\ge v$ in $B_R(z))$. Moreover, setting $r=|x-z|$, for $x\in B_R(z)$, we have
$$w(r)=v(x)\le m(|x-z|)=m(r),\;\;\mbox{and $\;h(m(r))\ge h(w(r)).$}$$
We also note that $v(z)\le u(z)$; hence, recalling from (\ref{h2h2}) that $w'\le 0$, we have that $\sup_{0\le r\le R} w(r)=w(0)=v(z)\le u(z).$
Thus from (\ref{h2h2}) and (\ref{vm}),
$$\Df w=-h(m(r))\le -h(w(r)),\;\;\mbox{$0<r<R,\;\;$ and $\;\;w(R)=\ell$}.$$
Multiplying both sides of the above equation by $w^{\prime}$ and integrating we have
$$ w^{\prime}(r)\le - \sqrt{2}\left(\int_{w(r)}^{w(0)}  h(s) ds\right)^{1/4}= -\sqrt[4]{4(H(w(0))-H(t))},$$
from which it follows that
$$\int_{w(r)}^{w(0)}\frac{dt}{\sqrt[4]{H(w(0))-H(t)}}\,dt\ge \sqrt{2}r.$$
Recalling (\ref{dd2}), we obtain the estimate
$$\sqrt{2}R\le \int_{\ell}^{w(0)}\frac{dt}{\sqrt[4]{H(w(0))-H(t)}}\,dt\leq \sup_{a>\ell}\zeta(a)=M_f<\infty.$$

It is now clear that if $R>M_f/\sqrt{2}$, then Problem (\ref{decrr}) can not have a solution $u\in C(\overline{\Om})$. $\;\;\;\Box$\\

\begin{rem}\label{eig2}
We now provide examples of functions $f\in C(\Om\times\R,[0,\infty))$ that satisfy (\ref{dd1}) and (\ref{dd2}) and hence the corresponding Problem (\ref{decrr}) has no solution in domains that are large enough. We look at two in particular. By Lemma \ref{app2} of the Appendix, and Theorem \ref{uss} above, we see that
if
$$\sup_{a>\ell}\frac{(a-\ell)^4}{H(a)}<\infty$$ then the Dirichlet Problem (\ref{decrr}) has no solution in domains whose in-radius is too large. We use this fact in the next two examples.  \\

\NI (i) As our first example we take $f(x,t)=h(t)=e^t$. By Lemma \ref{app2} of the Appendix, we note that
\begin{align*}
\int_{\ell}^a\frac{1}{\sqrt[4]{e^a-e^t}}\;dt
&=\zeta(a)=\psi(\ell)\\[.3cm]&\leq
\frac43\left(\frac{(a-\ell)^4}{H(a)}\right)^{1/4}=
\frac{4(a-\ell)}{3(e^a-e^{\ell})^{1/4}}\;\;\;\;\;\;\;\;a>\ell.
\end{align*}$\;$\\

It is clear that both (\ref{dd1}) and (\ref{dd2}) hold, and the Dirichlet problem (\ref{decrr}) has no solution if the in-radius $R$ of the domain $\Om$ is large enough.\\

\NI (ii) The next one we consider is the function
\eqRef{sf}
f(x,t)=\left\{\begin{array}{lcc} c t^{\al}+\dl & t\ge 0\\[.2cm] \dl & t<0 \end{array}\right.,
\end{equation}
where $c$ is a positive constant $\al\ge 3$ and $\dl>0$. Clearly, $h(t)=f(x,t)$ in this case. As an example, we consider boundary data with $\inf_{\p\Om} b=\ell\ge 0,$
and $b$ not identically zero on $\p\Om$. Clearly, solutions $u$ are positive in $\Om$, if they exist.
Here again, on using Lemma \ref{app2} of the Appendix we find that

$$\zeta(a)\le \frac{4(\al+1)^{1/4}(a-\ell)}{3[c(a^{\al+1}-\ell^{\al+1})+(\al+1)\dl(a-\ell)]^{1/4}},
\;\;\;\;\;\;a>\ell.$$
On noting that the right hand side tends to zero as $a\rar\ell^+$ and as $a\rar\infty$, we conclude that both (\ref{dd1}) and (\ref{dd2}) hold. Therefore, by Theorem \ref{uss} the Dirichlet problem (\ref{decrr}) has no solutions if the in-radius $R$ is large enough. The same conclusion holds if $\dt=0$ and $\ell>0$.  In the case $\dt=0$ and $\ell=0$, it is easily seen that
\eqRef{esss}
\zeta(a)\le \frac{4(\al+1)^{1/4}}{3c^{1/4}} a^{(3-\al)/4},\;\;\;a>0.
\end{equation}
Note that for $\al=3$, condition (\ref{dd2}) holds and therefore non-trivial solutions do not exist in $\Om$ when the in-radius $R$ is large enough. However if $\al>0$ and $\al\not=3$, the estimate (\ref{esss}) does not necessarily lead to condition (\ref{dd2}). Let us consider the situation $b\equiv 0$ on $\p\Om$, when the parameter $\dt$ in (\ref{sf}) is zero and when $\al>0$ with $\al\not=3$. Obviously $u\equiv 0$ is a solution to Problem (\ref{decrr}). Are there other solutions to (\ref{decrr})?  It turns out that a positive solution exists in $B_R(x_0)$ for any $R>0$ and $x_0\in\R^N$. See Remark \ref{alp} below. If $0<\al<3$, there are positive solutions of Problem (\ref{decrr}), regardless of the size of the domain. See Theorem \ref{apb2}, Remark 4.4 below, and Remark \ref{examp} in Section 6.\\
\end{rem}

\NI Next we prove a positive result in this context. We replace the condition (\ref{dd2}) by condition (\ref{dd3}) below. In a sense this complements the result in the preceding lemma. \\

\NI First we introduce some additional notations. Given $b\in C(\p\Om)$, we set $\ell:=\inf_{\p\Om}b$, and define
\eqRef{gg}
h_1(t):=\inf_{\Om\times [t,\;\infty)} f(x,s),\;\;\;\mbox{and}\;\;\;h_2(t)=\sup_{\Om\times[\ell,\;t]} f(x,s),\;\;
\mbox{for every $t\ge \ell$.}
\ee
We note that both $h_1$ and $h_2$ are non-decreasing, and $0\le h_1(t)\leq f(x,t)\leq h_2(t),\;\;t\ge \ell.$ We remark that if we set
$$\eta_1(t)=\frac{1}{t-\ell}\int_\ell^th_1(s)\,ds,\;\;\;\;\tx{and}\;\;\;\;
\eta_2(t)=\frac{1}{t-\ell}\int_t^{2t-\ell}h_2(s)\,ds,\;\;\;\;\;\;t>\ell$$
then $\eta_k$ are continuous in $(\ell,\infty)$ and a change of variable shows that each is non-decreasing in $(\ell,\infty)$. Moreover we see that $0\leq \eta_1(t)\leq h_1(t)$ and $\eta_2(t)\geq h_2(t)$ for $t>\ell$.

\noindent Therefore in what follows, we may assume that each $h_k$ is continuous in $(\ell,\infty)$, otherwise use $\eta_k$ in place of $h_k$.
For $a>\ell$, define
\eqRef{zta}
\zeta_k(a)=
\frac{1}{\sqrt{2}}\int_{\ell}^a\frac{1}{\sqrt[4]{H_k(a)-H_k(t)}}\,dt,
\;\;\;\;\;\;\;\tx{where}\;\;\;\;H_k(t)=\int_\ell^t h_k(s)\,ds\;\;\;\tx{for}\;\;k=1,2.
\ee$\;$\\

\NI We now require the following conditions on $f$:\\
\eqRef{dd3}\left\{
\begin{array}{ccl}
  (i)&\;&\liminf_{a\rightarrow \ell+}\zeta_1(a)=0  \;\;\;\;\;\;\;\tx{and} \\[.4cm]
 (ii)&\; &\limsup_{a\rightarrow \infty} \zeta_2(a)=\infty
\end{array}\right.
\ee

\nin In the next theorem we prove the existence of a solution $u\in C(\overline{\Om})$ to (\ref{decrr}) with $u>\ell$ in $\Om$. \\

\begin{lem}\label{apb2} Let $f:\Om\times\R\to[0,\infty)$ be a continuous function that satisfies condition (\ref{sup}). Let $h_k$ and $\zeta_k,\;k=1,2,$ be as defined in (\ref{gg}) and (\ref{zta}), respectively.
Assume further that conditions
(\ref{dd1}) and (\ref{dd3}) hold. Then Problem (\ref{decrr}) admits a solution $u\in C(\overline{\Om})$ with $u>\ell$ in $\Om$.

\end{lem}

\pf We will use Theorem \ref{sst} to prove the result. To this end we construct a sub-solution $w$ and a super-solution $v$ such that $w<v$ in $\Om$, with $w\ge \ell$ in $\Om$, $w\le b$, on $\p\Om$, and $v\ge \sup_{\p\Om}b$ on $\p\Om$.
\\

\nin First we proceed to construct the super-solution $v$, which will be a radial solution of the problem $\Df v=-h_2(v)$ in a ball that contains $\Om$. Fix $x^*\not\in\Om$. Because of (\ref{dd3}) we select $\bt>2\max_{\p\Om}b$, large enough such that $\Om\su B_{R^*/2}(x^*)$
where $R^*:=\zeta_2(\bt)$. We now consider the following decreasing function
$$\psi(t)=\frac{1}{\sqrt{2}}\int_t^\bt \frac{1}{\sqrt[4]{H_2(\bt)-H_2(s)}}\;ds\;\;\;\;\;\;\;\;\;\;\;\;\ell\leq t\leq \bt,$$
and let $\phi$ be the inverse of $\psi:[\ell,\bt]\to[0,R]$. Then by Lemma \ref{app4} in the Appendix, we note that
$$(\phi'(r))^2\phi''(r)=-h_2(\phi(r))),\;\;\;\;0<r<R^*.$$
Clearly, $\phi$ is decreasing, $\phi(0)=\bt$ and $\phi(r)\geq \ell,\;0\le r<R^*.$ We define
$$v(x):=\phi(|x-x^*|),\;\;\;\;\;\;\;\;\;x\in B_{R^*}(x^*).$$
It is obvious that $v\in C^2(\Om)$, and
$$\Df v=-h_2(v)\le -f(x, v)\;\;\;\tx{in}\;\;  \Om.$$
By observing that $\phi(r)$ is concave, it follows that
$$\phi(r)\ge \left(\frac{R^*-r}{R^*}\right)\bt,\;\;0\le r\le R^*.$$
Therefore $v(x)\geq \bt/2\geq \sup_{\p\Om}b$ for all $x\in\Om\su B_{R^*/2}(x^*)$.
We have thus shown that $v$ is a super-solution to (\ref{decrr}).\\

\NI Next we proceed to construct a sub-solution of (\ref{decrr}). If $R$ is the in-radius of $\Om$, by virtue of (\ref{dd3}) we choose
$$\ell<\al<\inf_\Om v\;\;\;\;\;\;\;\tx{such that}\;\;\;0<\zeta_1(\al)\leq R.$$
Let $R_*:=\zeta_1(\al)$. We note that $B_{R_*}(x_*)\su\Om$ for some $x_*\in\Om$.\\

\nin Here again, let us consider the decreasing function
$$\psi(t)=\frac{1}{\sqrt{2}}\int_t^\al\frac{1}{\sqrt[4]{H_1(\al)-H_1(s)}}\,ds
\;\;\;\;\;\;\;\ell\leq t\leq \al.$$
Let $\vp$ be the inverse of $\psi:[\ell,\al]\rar[0,R_*]$. As observed before, $\vp$ satisfies the equation
$$(\vp'(r))^2\vp''(r)=-h_1(\vp(r)),\;\;\;0<r<R_*,$$
and that $\ell\leq \vp(r)\leq \al$ in $[0,R_*]$ with $\vp(R_*)=\ell$, and $\vp(0)=\al$. We set
$$w(x)=\left\{\begin{array}{cl}
         \vp(|x-x_*|)&\;\;\;\;\;\;\;\;x\in B_{R_*}(x_*)\\[.2cm]
         \ell &  \;\;\;\;\;\;\;\;x\notin B_{R_*}(x_*).
       \end{array}\right.
$$
Note that as a result of our choice of $\al>\ell$ we see that
$$w<v\;\;\;\;\;\tx{in}\;\;\overline{\Om}.$$
We wish to show that $w$ satisfies $\Df w\ge -h_1(w)$ in all of $\R^N$, in the viscosity sense. So let $\psi\in C^2(\R^N)$ such that $w-\psi$ has a local maximum at $z\in\R^N$. It is enough to consider the case when either $z=x_*$, or $|z-x_*|=R_*$. We consider each case separately. Let us first consider the case when $z=x_*$. We notice that since $\vp'(0)=0$, we have $D\psi(x_*)=Dw(x_*)=0$. Therefore
$$\Df\psi(x_*)=0\geq -h_1(w(x_*))\geq -f(x_*,w(x_*)).$$
So we now assume that $|z-x_*|=R_*$. Note that in this case we have $w(z)=\vp(|z-x_*|)=\vp(R_*)=\ell$. We claim that $D\psi(z)=0$. Assume to the contrary, and let
$$e:=\frac{D\psi(z)}{|D\psi(z)|}.$$
First we note that, since $\ell\leq \vp(r)$, for $0\leq r\leq R_*$, and $w(z)=\ell$, we have
$$0\leq w(x)-w(z)\leq\psi(x)-\psi(z)=\left\langle D\psi(z),x-z\right\rangle+\frac12\left\langle D^2\psi(z)(x-z),x-z\right\rangle+o(|x-z|^2).$$
We now take $x:=z-\lam e$ for $\lam>0$ in the last inequality above. We obtain
$$\frac{\lam^2}{2}\left\langle D^2\psi(z)e,e\right\rangle+o(\lam^2)\geq \lam |D\psi(z)|.$$ Taking the limit as $\lam\rar0$ we get a contradiction. Therefore our claim holds and hence
$$\Df\psi(z)=0\geq -f(z,w(z)).$$
\nin Now we invoke Theorem \ref{sst} to conclude that Problem (\ref{decrr}) admits a solution $u\in C(\overline{\Om})$ such that $w\leq u\leq v$ in $\overline{\Om}$. Since $u(x_0)\geq w(x_0)=\al>\ell$, and $u$ is infinity super-harmonic in $\Om$, we conclude that $u>\ell$ in $\Om$.
\epf

\begin{rem} Here we would like to highlight the role that condition (\ref{dd3})(i) plays in ensuring that the solution $u$ in the above lemma
is strictly bigger than $\ell$ in a ball contained in $\Om$. Therefore, in the event that $b\equiv\ell$, and $f(x,\ell)\equiv0$ in $\Om$ for instance, the above lemma provides a non-trivial second solution to (\ref{decrr}).
\end{rem}

\begin{rem} Given $f\in C(\Om\times\R,[0,\infty))$ that satisfies condition (\ref{sup}), let $h_k,\;H_k$ and $\zeta_k$ be as defined in (\ref{gg}) and (\ref{zta}) above. If
\eqRef{gg1}
H_1^{-1/4}\in L^1(\ell\,,\,\ell+1),
\end{equation}
then $\zeta_1(a)\rar0$ as $a\to\ell^+$. Moreover, if
\eqRef{gg2}
h_2(t)=o(t^3)\;\;\;\;\;\tx{as $t\rar\infty$},
\end{equation} then $\zeta_2(a)\rar\infty$ as $a\rar\infty$. These assertions follow from Lemma \ref{app3} of the Appendix.
Therefore if $f$ satisfies conditions (\ref{gg1}) and (\ref{gg2}), then by Lemma \ref{apb2} above the Dirichlet problem (\ref{decrr}) admits a solution $u\in C(\overline{\Om})$.
Finally, for $f(x,t)=t|t|^{\g-1},\;0<\g<3,$ the function
$$v(x)=\left\{\frac{\s \left(R^{4/3}-|x|^{4/3}\right)}{\beta}\right\}^{\beta},\;\;0\le |x|\le R,$$
where $\beta=3/(3-\g)$, gives an explicit sub-solution to (\ref{decrr}) with $\Om:=B_R(o)$ and $b\equiv 0$. An extension by zero also provides a sub-solution in all of $\R^n$.
\end{rem}
\vsp

\section{\bf A priori $L^{\infty}$ Bounds and Existence of Solutions}

In this section we discuss a priori $L^{\infty}$ bounds for solutions to the Dirichlet problem in (\ref{dp}), that is, for solutions $u$ to
$$\Df u=f(x,u)\;\;\;\;\mbox{in $\Om$,\quad and\quad$u=b$ on $\p\Om.$}$$
Our main focus is establishing existence of solutions when the right-hand side function $f(x,t)$ satisfies some growth conditions. The main idea is that if the growth of $f(x,t)$ at infinity bears a certain relation to the function $|t|^3$, then we have a priori supremum bounds for solutions, see Theorem \ref{supf} below. These bounds are then shown to lead to the existence of solutions, see Corollary \ref{exs1}. In contrast to Theorem \ref{decf}, there are no restrictions on the underlying diameter of $\Om$. We point out that our result on a priori bounds, Theorem \ref{supf}, is optimal in the sense that if we choose $f(x,t)=-t |t|^{\dl},\;\dl\ge 2$, then a priori bounds may no longer be available, see Lemma \ref{apb1}, in Section 6. However, solutions may still continue to exist. The discussion in this section leads to a fairly general class of functions $f$ for which Problem (\ref{dp}) admits a solution. In a sense, the results of this section complement those obtained in Section 3. \\

We address first the case when $f(x,t)=h(x)$.  That is we consider the Dirichlet problem
\eqRef{h}
\left\{\begin{array}{rcll}
         \Df u &=&h(x)&\;\;x\in\Om \\[.3cm]
         u &=&b&\;\;x\in\p\Om.
       \end{array}
 \right.
\end{equation}

For any function $g$, define $g^+=\max\{g,0\}$ and $g^-=\min\{g,0\}$.\\

\begin{thm}\label{bounds1}
Let $\Om\su \R^N$ be a bounded domain whose out-radius is R . Let $h\in C(\Om)\cap L^\infty(\Om)$, and $b\in C(\p\Om)$.
\begin{enumerate}[(i)]
\item If $u\in C(\overline{\Om})$ is a sub-solution of (\ref{h}), then
$$u\le \sup_{\partial \Om}b-\s(\inf_{\Om} h^-)^{1/3}R^{4/3}\;\;\;\;\tx{in}\;\;\overline{\Om}.$$
\item If $u\in C(\overline{\Om})$ is a super-solution of (\ref{h}), then
$$u\geq  \inf_{\partial \Om}b-\s(\sup_{\Om} h^+)^{1/3}{R}^{4/3}\;\;\;\;\tx{in}\;\;\overline{\Om}.$$
\end{enumerate}
\end{thm}

\pf Let $u\in C(\overline{\Om})$ be a sub-solution of (\ref{h}). Fix $\ep>0$, and suppose that $\Om\su B_R(z_0)$ for some $z_0\in\Om$. We consider the function
$$w(x)= \sup_{\partial \Om}b+\s(-\inf h^- +\ep)^{1/3}\left(R^{4/3}-|x-z_0|^{4/3}\right),\;x\in B_{R}(z_0)$$
It is easy to show that $\Df w=\inf h^- -\ep<h$ in $\Om$, and since $\partial \Om\subset \overline{B}_{R}(z_0)$, $u\le \sup_{\partial \Om}b\le w$ on $\partial \Om$.
By the Comparison Principle, Lemma \ref{mcomp}, we get
\be\label{p}
u(x)\le  \sup_{\partial \Om}b+\s(-\inf h^- +\ep)^{1/3}\left(R^{4/3}-|x-z_0|^{4/3}\right),\;x\in \Om.
\ee
Since $\ep>0$ is arbitrary, this proves (i). To get the estimate in (ii), let $u\in C(\overline{\Om})$ be a super-solution of (\ref{h}). Then $-u$ is a sub-solution of $\Df v=-h$ in $\Om$ with $v=-b$ on $\p\Om$. Thus the estimate in (ii) is obtained by an application of the estimate in (i) to $-u$.  \qed\\

We have the following immediate corollary.

\begin{cor}\label{ccc}Let $\Om\su \R^N$ be a bounded domain whose out-radius is $R$. Suppose $h\in C(\Om)\cap L^\infty(\Om)$, and $b\in C(\p\Om)$. For any solution $u\in C(\Om)$ of (\ref{h}) the following estimate holds:
$$\inf_{\partial \Om}b-\s(\sup_{\Om} h^+)^{1/3}{R}^{4/3}\leq  u \le \sup_{\partial \Om}b-\s(\inf_{\Om} h^-)^{1/3}{R}^{4/3}\;\;\;\;\tx{in}\;\;
\overline{\Om}.$$
\end{cor}$\;$\\

We now extend the above a priori estimate to solutions of the Dirichlet problem (\ref{dp}) when the inhomogeneous term $f(x,t)$ is subjected to appropriate growth conditions as $|t|\rar\infty$.\\

\begin{thm}\label{supf}
 Let $\Om\su\R^N$ be a bounded domain, $b\in C(\p\Om)$ and $f\in C(\Om\times\R,\R)$ such that (\ref{sup}) holds. Assume that
\eqRef{cnn}
\left\{\begin{array}{ll}
         \displaystyle{(\tx{i})} & \displaystyle{\liminf_{t\rar\infty}\frac{\inf_\Om f(x,t)}{t^3}:=\bt} \\[.4cm]
          \displaystyle{(\tx{ii})} &\displaystyle{ \liminf_{t\rar-\infty}\frac{\sup_\Om f(x,t)}{t^3}:=\al}
       \end{array}\right.
\end{equation}
for some $\al,\bt\in[0,\infty]$. Then there is a positive constant $C$, depending on $f, b$ and diam$(\Om)$ such that $\|u\|_{L^{\infty}(\Om)}\leq C$ for any solution $u\in C(\overline{\Om})$ of (\ref{dp}).
\end{thm}

\pf
\noindent We start with condition (\ref{cnn})(i) to establish an upper bound. Let $\ep>0$, to be chosen later. Then according to (\ref{cnn})(i) there is $t_\bt>0$,  depending on $f$ and $\ep$, such that
\eqRef{ine}
f(x,t)\geq -\ep t^3,\;\;\;(x,t)\in\Om\times[t_\bt,\infty).
\end{equation}
Given a viscosity solution $u\in C(\overline{\Om})$ of (\ref{dp}), we consider the following open subset of $\Om$.
$$\Om_\bt:=\{x\in\Om:u(x)>t_\bt\}.$$
We assume that $\Om_\bt$ is non-empty. Without loss of generality we may suppose $t_\bt\geq \sup_{\p\Om}b.$
Clearly $t_\bt< M$, where $M:=\sup_\Om u$. Let $v\in C(\overline{\Om})$ be the solution of the following Dirichlet problem.

\eqRef{sup4}
\Df v=-\ep M^3\;\;\mbox{in} \;\;\Om,\;\;\mbox{and}\;\;\;v=t_\bt\;\;\mbox{on}\;\;\p \Om.
\end{equation}
By the Comparison Principle, Lemma \ref{luw}, we note that $v\geq t_\bt$ in $\Om$. By Theorem \ref{bounds1}, we see that
\eqRef{ves}
\sup_\Om v\leq t_\bt+\ep^{1/3}\s MR^{4/3},
\end{equation}
where $R$ is the out-radius of $\Om$.
Now, note that  in the sub-domain $\Om_\bt$ we have
$\Df u=f(x,u)\geq -\ep u^3\geq-\ep M^3$ with $u=t_\bt$ on the boundary $\p\Om_\bt$. Therefore, by Lemma \ref{luw}, we find that $u\leq v$ in $\Om_\bt$, and thus $M\leq \sup_\Om v$.\\

\noindent Consequently this last estimate together with (\ref{ves}) gives
\eqRef{wv}
M\leq \sup_\Om v\leq t_\bt+\ep^{1/3}\s M R^{4/3}.
\end{equation}
We choose $\ep=(2\s R^{4/3})^{-3}$, and the inequality in (\ref{wv}) then leads to

\eqRef{gl} u\leq  M\leq 2t_\bt,
\end{equation}
provided $\Om_\bt$ is non-empty. Since $u\leq t_\bt$ if $\Om_\bt=\emptyset$, we conclude that the estimate (\ref{gl}) actually holds in any case.\\

\nin Next we use (\ref{cnn})(ii) to get a lower bound. As before, corresponding to $\ep:=(2\s R^{4/3})^{-3}$ we can find $t_\al<\min\{0,\inf_{\p\Om}b\}$ such that
$$f(x,t)\leq -\ep t^3\;\;\;\;\;\;(x,t)\in\Om\times(-\infty,t_\al].$$
Let us observe that $w:=-u$ is a solution of
$$\Df w=g(x,w)\;\;\;\mbox{in}\;\;\Om\;\;\;\mbox{and}\;\;\;w=-b\;\;\;
\mbox{on}\;\;\p\Om,
$$ where $g(x,t)=-f(x,-t)$ satisfies
$$g(x,t)\geq -\ep t^3\;\;\;\;\;\;\;(x,t)\in\Om\times[-t_\al,\infty).$$  Therefore by the above argument we find that
$$-u=w\leq -2t_\al.$$  Therefore $u\geq 2t_\al$ on $\Om$. In summary, we find that there are constants $t_\al<\min\{\inf_{\p\Om} b,0\}\leq \max\{0,\sup_{\p\Om}b\}<t_\bt$ such that
$$2t_\al\leq u(x)\leq 2t_\bt,\;\;\;\;\;\;(x\in\overline{\Om}).$$
 \epf\\

\begin{rem} The class of functions $f(x,t)$ described in (\ref{cnn}) include those $f:\Om\times\R\rar\R$ that
\begin{enumerate}[(i)]
\item are non-decreasing in $t$ for each $x\in\Om$, or

\item satisfy Condition (\ref{sup}) and for some $t_{-}<t_{+}$
$$f(x,t)\leq 0\;\;\;\forall\;(x,t)\in\Om\times(-\infty,t_{-})\;\;\;\;\;
\tx{and}\;\;\;\;\;f(x,t)\geq 0\;\;\;\forall\;(x,t)\in\Om\times(t_{+},\infty),\;\;\;\tx{or}$$
\item satisfy Condition (\ref{sup}) and that satisfy $f(x,t_n)=0=f(x,s_n)$ for some sequence $\{s_n\}$ and $\{t_n\}$ with $t_n\rar\infty$ and $s_n\rar-\infty$. For instance $f(x,t)=a(x)(1+\cos t)^\g$ for $\g>0$ and $a\in C(\Om)\cap L^\infty(\Om)$.
\end{enumerate}
\end{rem}$\;$\\

\noindent We now state our main existence result of this section.\\

\begin{thm}\label{exs1}
Let $\Om\su \R^N$ be a bounded domain, and $b\in C(\partial \Om)$. If $f\in C(\Om\times\R,\R)$ satisfies the conditions in (\ref{sup}) and (\ref{cnn}), then the Dirichlet problem (\ref{dp}) admits a solution $u\in C(\overline{\Om})$.

\end{thm}
\vsp
\pf Let $R$ be the out-radius of $\Om$ and $\ep=(2\s R^{4/3})^{-3}$. Since $f$ satisfies (\ref{cnn}), we can find constants $t_\al<\min\{0,\inf_\Om b\}$ and $t_\bt> \max\{0,\sup_{\p\Om}b\}$ such that
\eqRef{fin}
f(x,t)\leq -\ep t^3,\;\;\;(x,t)\in\Om\times(-\infty,t_\al]\;\;\;\;\tx{and}
\;\;\;\;\;
f(x,t)\geq -\ep t^3,\;\;\;(x,t)\in\Om\times[t_\bt,\infty).
\end{equation}
We take $C\geq \max\{-2t_\al,2t_\bt\}$ and we now define
$$\widehat{f}(x,t)=\left\{\begin{array}{lll}
f(x,C)&\;\;\;\tx{if} & \;t>C\\[.2cm]f(x,t) &\;\;\;\tx{if} &\;|t|\le C\\[.2cm]f(x,-C) &\;\;\;\tx{if} & \;t<-C.\end{array}\right.$$
Clearly, as a result of condition (\ref{sup}), we see that $\widehat{f}\in L^{\infty}(\Om\times\R)$. Therefore Corollary \ref{bd} applies, implying the existence of a solution $u\in C(\overline{\Om})$ to the problem
$$\Df u=\widehat{f}(x,u)\;\;\mbox{in $\Om$ and $u=b$ on $\partial \Om$.}$$

The choice of the constant $C$ shows that the inequality (\ref{fin}) holds with $f$ replaced by $\widehat{f}$, but with the same constants $t_\al$ and $t_\bt$. Therefore, as in the proof of Theorem \ref{supf} above, we conclude that
$$2t_\al\leq u(x)\leq 2t_\bt\;\;\;\;\;\;\;(x\in\Om).$$
But this implies that $|u(x)|\leq C$ in $\Om$, and therefore $u$ is a solution of (\ref{dp}) as desired.\qed $\\$

As pointed out in Section 1, we have been applying Theorem \ref{sst} to various situations in order to conclude the existence of a solution. So far we have not discussed the question as to whether multiple solutions may exist in some instances. Clearly this is intimately connected to the question of uniqueness and the examples we discuss in Section 6 indicate lack of uniqueness in general. In this connection we have included the following remark, specifically, to indicate that multiplicity of solutions is not an unusual occurrence. \\

\begin{rem}\label{eig}
In some instances, Lemma \ref{apb2} leads to the existence of a second solution which is not guaranteed by Theorem \ref{exs1}.
As an example, consider the Dirichlet problem with a non-positive right hand side, that is,
$$\Df u=-f(x,u)\;\;\mbox{in $\Om$ with $u=c$ on $\p\Om$}, $$
where $c$ is a constant and $f\in C(\Om\times\R,[0,\infty))$. Suppose that $f(x,c)=0$ for all $x\in \Om$. Then clearly one solution is $u_1=c$ in $\Om$. Suppose that $f$ satisfies the hypotheses of Lemma \ref{apb2}. The sub-solution $w$ constructed in Lemma \ref{apb2} satisfies $w>c$ in a ball $B\su\Om$, implying that Lemma \ref{apb2} provides us with a second solution $u_2>c$ in a ball $B\su\Om$. Thus we obtain two solutions. Also see Remark \ref{examp} in this context.\\

\end{rem}

\section{\bf An Example on lack of a priori bounds}
\vsp

In this section we look at the following Dirichlet problem on a bounded domain $\Om$.
\eqRef{disc}
\lf\{\begin{array}{rcll}
         \Df u &=& f(x,u) &\;\mbox{in $\;\Om$}\\[.2cm]
         u&= &0  &\;\mbox{on $\;\p \Om$},
       \end{array}
 \right.
\end{equation}$\;$\\

Our primary objective is to demonstrate that a priori bounds of the kind  announced in Theorem \ref{supf} may not hold if the hypotheses stated in (\ref{cnn}) are not met. In fact, if we take $f(x,t)=-t|t|^{\g-1},\;\g>3$, Theorem \ref{supf} does not apply. As a matter of fact, as we will see below, a priori $L^{\infty}$ bounds do not exist in these cases. This will be shown by first considering Problem (\ref{disc}) in balls.

Finally we include a brief discussion of Problem (\ref{disc}) when $f(x,t)=-a(x)t^3$ where $a\in C(\Om)\cap L^\infty(\Om)$.\\

Consider the following Dirichlet problem in the unit ball $B=B_1(o)$, centered at the origin,

\eqRef{dpb}
\lf\{\begin{array}{rcll}
         \Df u &=& -u|u|^{\g-1} &\;\mbox{in $\;B$}\\[.2cm]
         u&= &0  &\;\mbox{on $\;\p B$},
       \end{array}
 \right.
\end{equation}

\begin{lem}\label{apb1} Let $\g>3$, then Problem (\ref{dpb}) admits a sequence $\{u_k\}$ of solutions in $C(\overline{B})$ such that
$$\lim_{k\rar\infty}\|u_k\|_{L^\infty(B)}=\infty.$$
\end{lem}

\pf To prove the lemma, we will construct radial solutions $u_k$ with increasing $L^{\infty}$ norm.

We first carry out some elementary computations that will be useful in the construction of this sequence of solutions. A change of variables shows that, for $a>0$,
\eqRef{est}
\int_{0}^{a}\frac{ds}{(a^{\g+1}-s^{\g+1})^{1/4}}=
\frac{1}{a^{(\g-3)/4}}\int_0^1\frac{dt}{(1-t^{\g+1})^{1/4}}.
\end{equation}
Since the right hand side exists, one can choose a unique $a=a(\g)>0$ such that
\eqRef{aaa}\int_{0}^a\frac{ds}{(a^{\g+1}-s^{\g+1})^{1/4}}=
\left(\frac{4}{\g+1}\right)^{1/4}.
\end{equation}
Hereafter we take $a$ to be the unique value $a(\g)$ determined by the above equality.
Let $\vp$ be the inverse of the decreasing function $\psi:[0,a]\rar[0,1]$ given by
$$\psi(t)=\left(\frac{\g+1}{4}
\right)^{1/4}\int_t^a\frac{ds}
{(a^{\g+1}-s^{\g+1})^{1/4}}$$
For $0\le r\le 1$, we have
\eqRef{sol}
\int_{\vp(r)}^a\frac{ds}
{(a^{\g+1}-s^{\g+1})^{1/4}}=\left(\frac{4}{\g+1}
\right)^{1/4}r.
\end{equation}$\;$\\

\nin We observe that for $0<r<1$,
\be\label{sun}
\left\{\begin{array}{lcr}\vp'(r)=-
\left(\frac{4}{\g+1}\right)^{1/4}\left(a^{\g+1}-\vp(r)^{\g+1}
\right)^{1/4},\\[.3cm]
\vp''(r)=-\left(\frac{\g+1}{4}\right)^{1/2} \vp(r)^\g\left(a^{\g+1}-\vp(r)^{\g+1}\right)^{-1/2}. \end{array} \right.
\ee
Thus $\vp\in C^2((0,1])\cap C^1([0,1])$ is decreasing and concave on $[0,1]$, with $$\vp'(0+)=0,\;\;\vp(1)=0,\;\;\vp^{'}(1-)<0\;\;\;\tx{and}\;\;\; \vp''(1-)=0.$$ Moreover we see that $\vp$ is a solution of
the initial value problem
$$
\lf\{\begin{array}{rcl}
         (\vp'(r))^2\vp''(r) &=& -\vp|\vp|^{\g-1}(r), \;\;\;\;0<r<1\\[.3cm]
         \vp(0)&= &a, \; \;\vp'(0)=0,\;\;\vp(1)=0.
       \end{array}
 \right.
$$

\nin We now extend $\vp$ to the interval $[0,4]$ as follows.
$$\vp_*(r):=\left\{\begin{array}{rcl}
                   \vp(r) & &0\leq r\leq 1 \\[.2cm]
                   -\vp(2-r) & &1\leq r\leq 2 \\[.2cm]
                   -\vp(r-2)&&2\leq r\leq 3\\[.2cm]
                   \vp(4-r)&&3\leq r\leq 4.
                 \end{array}
 \right. $$
\nin Thus $\vp_*$ is obtained by first carrying out an odd extension of $\vp$ to $[0,2]$, and then carrying out an even extension of the resulting function on $[0,2]$, to $[0,4]$. Before we proceed further, let us make two simple observations about the extension $\vp_*$. First we note that $\vp_*$ satisfies the equation $(\vp_*')^2\vp_*''=-\vp|\vp|_*^{\g-1}$ in the set $S=(0,2)\cup(2,4)$. This follows since $\vp_*$ is $C^2$ on this set, and $\vp_*(1)=\vp_*''(1)=0=\vp_*''(3)=\vp_*(3)$. Furthermore, we also note that $\vp_*'(0)=\vp_*'(2)=\vp_*'(4)=0$. We now extend $\vp_*$ to the set $[0,\infty)$ as a $4$-periodic function $\vp_\infty$. Hence, $\vp_{\infty}$ satisfies $(\vp_{\infty}')^2\vp_{\infty}''=-\vp_{\infty}|\vp_{\infty}|^{\g-1}$ in the set $\cup_{k=0}^{\infty}(2k,\;2k+2).$

Let us now define $w$ on $\R^N$ by
$$w(x):=\vp_\infty(|x|),\;\;\;\;\;\;\;\;x\in\R^N.$$ Note that $w\in C^1(\R^N)$.
First we show that $w$ is a viscosity solution of $\Df w=-w|w|^{\g-1}$ in $\R^N$. Let $\psi\in C^2(\R^N)$  such that $w-\psi$ has a local extrema at $x_0\in\R^N$. Since $\vp_{\infty}$ is $4$-periodic, and in view of the remarks we have made above, it is enough to assume that $x_0=o$ or $|x_0|=1$ or $|x_0|=2$.

In what follows, $e_i$ will denote the unit vector along the positive $x_i$ axis, for $i=1,2,\cdots, N$. We first consider the cases when $x_0=o$ and when $|x_0|=2$.
Recalling that $\vp_*'(0)=\vp_*'(2)=0$ we have $D w(x_0)=0$. Let us first take up the case $x_0=o$. So suppose that $w-\psi$ has a maximum at $x_0=o$.
Then $D\psi(o)=D w(o)=o$. Since $w(o)=\vp_{\infty}(0)=a>0$, we have that $0=\Df\psi(o)\ge -w(o)^\g$.  Now let $w-\psi$ have a minimum at $o$. Once again $D\psi(o)=o$ and
\begin{align*}
\vp_\infty(|x|)-\vp_\infty(0)&=w(x)-w(o)\\&\ge \psi(x)-\psi(o)=\frac12\langle D^2\psi(o)x,\;x \rangle+\circ(|x|^2).
\end{align*}
We take $x=\vep e_i$ for a sufficiently small $\vep>0$. Then
\eqRef{de}
D_{ii}\psi(o)\leq 2\left(\lim_{\vep\rightarrow 0}\frac{\vp_\infty(\vep)-\vp_\infty(0)}{\vep^2}\right)=2\left(
\lim_{\vep\rightarrow 0}\frac{\vp(\vep)-\vp(0)}{\vep^2}\right).
\end{equation}
By concavity, $(a^{\g+1}-t^{\g+1})^{1/4}\geq a^{\g/4}(a-t)^{1/4}$, for $0\leq t\leq a$, and hence for $0<r<1$, we estimate
$$
\int_{\vp(r)}^a\frac{dt}{(a^{\g+1}-t^{\g+1})^{1/4}}
\leq
\frac{1}{a^{\g/4}}\int_{\vp(r)}^a
\frac{dt}{(a-t)^{1/4}}
=\frac{4(a-\vp(r))^{3/4}}{3a^{\g/4}}.$$
Using this inequality in (\ref{sol}) and noting that $\vp(0)=a$, we find that
\eqRef{aa}\vp(0)-\vp(r)\geq \sigma\left(\frac{a^\g}{\g+1}\right)^{1/3}r^{4/3},\;\;\;\;\;\;
\text{where}\;\;\s=\frac{3\sqrt[3]{3}}{4}.
\end{equation}
We use (\ref{aa}) in (\ref{de}) to conclude that
$$D_{ii}\psi(o)\le \lim_{\vep\rightarrow 0}\frac{\vp(\vep)-\vp(0)}{\vep^2}
\leq -\sigma
\lim_{\vep\to0}\left(\frac{a^\g}{\vep^2(\g+1)}\right)^{1/3}=-\infty.
$$
This contradicts that $\psi$ belongs to $C^2$, and thus $x_0=o$ can not be a point of minimum of $w-\psi$.

Next we take up the case $|x_0|=2$. Since the arguments are analogous, we will be brief. Suppose then that $w-\psi$ has a local minimum at $x_0$. Then as before we note that
$$0=\Df\psi(x_0)\leq |\vp(0)|^\g=-\vp_*(2)|\vp_*(2)|^{\g-1}=-\vp_\infty(2)|\vp_{\infty}(2)|^{\g-1}=-w(x_0)|w(x_0)|^{\g-1}.$$
Now suppose that $w-\psi$ has a local maximum at $x_0$. Then, since $D\psi(x_0)=0$, as before
\begin{align*}
\vp_\infty(|x|)-\vp_\infty(2)&=w(x)-w(x_0)\\&\le \psi(x)-\psi(x_0)=\frac12\langle D^2\psi(x_0)(x-x_0),\;(x-x_0) \rangle+\circ(|x-x_0|^2).
\end{align*}
We now take $x=(1-\vep/2)x_0$ for a sufficiently small $\vep>0$. Then we have
$$\frac18\langle D^2\psi(x_0)x_0,\;x_0 \rangle\geq \lim_{\vep\to0}\frac{\vp_\infty(2-\ep)-\vp_\infty(2)}{\vep^2}=
\lim_{\vep\to0}\frac{\vp(0)-\vp(\vep)}{\vep^2}$$
which again, in view of (\ref{aa}), leads to a contradiction.\\

Finally, we take up the case when $|x_0|=1$. Let us recall that $\vp_{\infty}(1)=\vp_{\infty}''(1)=0$, and that $w\in C^2(A)$ where $A:=\{x\in\R^N:0<|x|<2\}$. For $x\in A$, we compute
\begin{align*}
\Df w(x)&=(\vp_\infty'(|x|))^2\vp_\infty''(|x|)|D (|x|)|^4+(\vp_\infty'(|x|)^3\Df(|x|)\\&
=(\vp_\infty'(|x|))^2\vp_\infty''(|x|).
\end{align*}
Thus evaluating $\Df w$ at $x_0$, we find that
$$\Df w(x_0)=(\vp_\infty(1))^2\vp_\infty''(1)=0.$$
Therefore $\Df w(x_0)=-w(x_0)|w(x_0)|^{\g-1}$.
In summary, we have shown that $w$ is a viscosity solution of $\Df w=-w|w|^{\g-1}$ in $\R^N$.\\

Finally we define a sequence of viscosity solutions of $\Df u=-u|u|^{\g-1}$ in $B$ such that $u=0$ on the boundary $\p B$. For each positive integer $k$ we set
$$u_k(x)=(2k-1)^{4/(\g-3)}w((2k-1)x),\;\;\;\;\;\;x\in B.$$
It is clear that $u_k$ is a viscosity solution of (\ref{dpb}) such that $u_k(o)=a(2k-1)^{4/(\g-3)}$, as desired. \epf
\\

\begin{rem}\label{alp}
Let $R>0$ be given, and suppose $\g>0,\;\g\not=3$. Choose $a=a(\g,R)$ such that the left-hand side of (\ref{aaa}) is equal to
$$\left(\frac{4}{\g+1}\right)^{1/4}R.$$
It is easily seen, by employing a change of variables, that $a^{\g -3}R^4=k$ where $k=k(\g)$ depends only on $\g$.

Let $\vp_R$ be defined on $[0,R]$ by
$$ \int_{\vp_R(r)}^a\frac{ds}
{(a^{\g+1}-s^{\g+1})^{1/4}}=\left(\frac{4}{\g+1}
\right)^{1/4}r.$$
Then as in the proof of Lemma \ref{apb1} one can verify that $v(x)=\vp_R(|x-x_0|)$ is a positive solution of
$$\Df v=-v|v|^{\g-1}\;\;\;\tx{in}\;\;B_R(x_0)\;\;\;\tx{and}\;\;v=0\;\;
\tx{on}\;\;\p B_R(x_0).$$
\end{rem}$\;$\\

\begin{rem} In contrast, recall part (ii) of Remark \ref{eig2}, in this context, as solutions, constructed in Lemma \ref{apb1}, change sign in $B$.   \\
\end{rem}

In the following remark, we discuss the case when the right hand side of (\ref{dpb}) does not grow any faster than the power three.

\begin{rem}\label{examp} We consider the following two cases.\\

\NI {\bf Case 1}: We take $f(x,t)=o(t^3)$, as $|t|\rightarrow \infty$, and for every $x\in \Om,\;f(x,0)=0.$ We now recall the discussion in Remark \ref{eig}. Clearly, $u_1(x)=0$ is a solution, and, since the condition (\ref{dd3}) is satisfied, a second solution $u_2(x)>0$ can be found by applying Lemma \ref{apb2}. Thus there are at least two solutions in $\Om$ to (\ref{disc}). These conclusions hold regardless of the sign of $f$.\\

\NI {\bf Case 2}: Consider the Dirichlet problem
\eqRef{eig1}
\Df u=-a(x)u^3\;\;\;\mbox{in $\Om\;$ with $u=b$ on $\p\Om$},
\ee
where $a(x)\in C(\Om)\cap L^{\infty}(\Om)$ and $b\in C(\p\Om)$. If $u\in C(\overline{\Om})$ is a solution of (\ref{eig1}) then by the estimates
in Theorem \ref{bounds1}, we have
$$
\inf_{\p\Om}b-\sigma (\sup_{\Om}|a(x)|M^3)^{1/3} R^{4/3}\le -M<M\le \sup_{\p\Om}b +\sigma (\sup_{\Om}|a(x)|M^3)^{1/3} R^{4/3}.
$$
where $M=\sup_{\Om}u$ and $R$ is the the out-radius of $\Om$. It follows that if $\sigma (\sup_{\Om}|a(x)|)^{1/3} R^{4/3}<1,$ then
\eqRef{infsup}
M\le \frac{\max(-\inf_{\p\Om}b,\;\sup_{\p\Om}b)}{1-\sigma (\sup_{\Om}|a(x)|)^{1/3} R^{4/3}}<\infty.
\end{equation}
In the event $b\equiv 0$ on $\p\Om$, we obtain that $u\equiv 0$ is the only solution, under the condition $\sigma(\sup_{\Om}|a(x)|)^{1/3} R^{4/3}<1$. This is a result that is quite general in nature and applies irrespective of the sign of the function $a(x)$. This observation is relevant only when $a(x)$ changes sign or when $a(x)\ge 0$. Notice that if $a\le 0$, in $\Om$, then $u$ is identically zero regardless of the size of the domain. This follows since $u$ would be infinity sub-harmonic wherever $u(x)>0$ and infinity super-harmonic wherever $u(x)<0$. \\

\NI A question that arises in this context is, under what conditions on $a(x)$ will (\ref{eig1}) have a nontrivial solution when $b\equiv 0$ on $\p\Om$. Consider the related problem
\eqRef{ju}
\Df u=-\lambda a(x) u^3,\;\;\mbox{in $\Om$ with $u=0$ on $\p\Om$},
\end{equation}
where $\lm\ge 0$ is a parameter, and $a(x)\geq 0$ and $a(x)\not\equiv 0$ in $\Om$. As was pointed out above, from (\ref{infsup}) we conclude that $u\equiv 0$ is the only solution of (\ref{ju}) in $C(\overline{\Om})$ if $\lm \sigma^3 \sup_{\Om}a(x) R^4<1$. Recently we became aware of the work \cite{JUT} in which P. Juutinen considers the eigenvalue problem $\Df^N u=-\lam u$ in bounded domains, where $\Df^Nu:=|Du|^{-2}\Df u$ is the so-called normalized infinity Laplacian. The methods of \cite{JUT} can easily be adapted to the eigenvalue Problem (\ref{ju}) with $a(x)\equiv1$. In Part II of the Appendix, we make a remark on the eigenvalue problem (\ref{ju}) when $a(x)\geq 0$ and $a(x)\not\equiv0$. See Proposition \ref{egn1}. It would be interesting to study the eigenvalue problem (\ref{ju}) when $a(x)$ is a non-negative and a non-trivial continuous function. \\
\\
\end{rem}

\section{\bf The Harnack Inequality and a Comparison Principle}

Our main goal in this section is to prove two results that may be of some independent interest. The first is a version of the Harnack inequality that applies to inhomogeneous equations. We also include a version of the comparison principle when $f(x,t)=f(t)$ and is non-decreasing, continuous and has no sign restrictions. The statement is a partial result of what one would hope to be the full comparison principle under the only assumptions that $f(t)$ is continuous and non-decreasing.

We start with the following Harnack inequality that provides a generalization of the well known Harnack inequality for non-negative infinity super-harmonic functions, see \cite{AMJ,BHA,CRM,LIM}.\\

\begin{thm} Suppose $h\in C(\Om)\cap L^\infty(\Om)$, and $u\in C(\Om)$ is a non-negative viscosity super-solution of $\Df u=h(x)$ in $\Om$. Then for any $z\in \Om$ such that  $B_{2r}(z)\su\Om$ we have
\eqRef{har}
\sup_Bu\leq 9\inf_B u+12\s\left(r^4\sup_\Om h^+\right)^{1/3}
\end{equation}
where $B:=B_{2r/3}(z)$.

\end{thm}
\pf Let $u\in C(\Om)$ be a positive function that satisfies $\Df u\leq h(x)$ in $\Om$ in the viscosity sense.\\

Consider a ball $B_{2r}(z)\su\Om$. Define $w(x):=\vp(|x-z|),\;x\in B_r(z),$ where
$$\vp(t):=u(z)\left(\frac{r-t}{r}\right)+\s(\sup h^+ )^{1/3}\left(t^{4/3}-r^{4/3}\right),\;\;\;\;t\in\R.$$
First we will show that
\eqRef{temp}
u(x)\geq \frac{1}{3}u(z)-\s(r^4\sup_\Om h^+)^{1/3},\;\;\;x\in B_r(z).
\end{equation}
For this we assume that
\eqRef{asu}
u(z)\geq3\s(r^4\sup_\Om h^+)^{1/3},
\end{equation}
for otherwise the inequality in (\ref{temp}) is trivially true. Since $\vp$ is radial, it is easily checked that, for $x\in B_r(z)\setminus\{z\}$ and $t:=|x-z|$,
\ben
\Df w&=\vp''(t)(\vp'(t))^2.
\een
Employing (\ref{asu}) we see that
\begin{align*}
(\vp'(t))^2\vp''(t)&=\frac49\s(t^{-2}\sup_\Om h^+)^{1/3}\left[\left(\frac{u(z)}{r}\right)^2-
\frac83\frac{u(z)}{r}\s(t\sup_\Om h^+)^{1/3}+\frac{16}{9}\s^2(t\sup_\Om h^+)^{2/3} \right]\\&
\geq \sup_\Om h^+.
\end{align*}
This, together with (\ref{asu}), shows that
\eqRef{su}
\Df w(x)\geq \sup_\Om h^+,\;\mbox{ for $x\in B_r(z)\setminus\{z\}$, \;\;$w(z)\le u(z)$ and $w=0$ on $\p B_r(z)$}.
\end{equation}

On the other hand $\Df u\leq h(x)\leq \sup_\Om h^+$ in $B_r(z)$ and since $u\geq 0$ in $\Om$ we see that
$u\geq w$ on the boundary of the punctured ball $B_r(z)\setminus \{z\}$. By the comparison principle in Lemma \ref{luw}, we see that
$w\leq u $ in $B_r(z)$, that is
\eqRef{bi}
u(x)\geq u(z)\left(\frac{r-|x-z|}{r}\right)+\s(\sup h^+ )^{1/3}\left(|x-z|^{4/3}-r^{4/3}\right),\;\;\;\;\;x\in B_r(x).
\end{equation}
Thus, for $x\in B_{2r/3}(z)$ we have
$$u(x)\geq \frac{u(z)}{3}-\s(r^4\sup_\Om h^+)^{1/3}.$$ In particular,
$$u(x_0)\geq \frac{u(z)}{3}-\s(r^4\sup_\Om h^+)^{1/3},$$
where $u(x_0)=\min_{B_{2r/3}(z)} u$.
Now let $x\in B_{2r/3}(z)$, so that $B_{2r/3}(z)\su B_{4r/3}(x)$. Applying (\ref{bi}) in the ball $B_{4r/3}(x)$, we observe that
\begin{align*}
u(x_0)&\geq\frac{u(z)}{3}-\s(r^4\sup h^+)^{1/3}
\geq \frac13\left[\frac{u(x)}{3}-\s(r^4\sup h^+)^{1/3}\right]-\s(r^4\sup h^+)^{1/3}\\
&=\frac19u(x)-\frac43\s(r^4\sup h^+)^{1/3}
\end{align*}
Therefore
$$u(x)\leq 9u(x_0)+12\s(r^4\sup_\Om h^+)^{1/3},\;\;\;x\in B_{2r/3}(z).$$ This concludes the proof of (\ref{har}).\epf\\

\NI{\bf Comparison}\\

\noindent For the next comparison theorem, let $f\in C(\R,\R)$ be a non-decreasing function. We define two extended real numbers $\ell_*$ and $\ell^*$ as follows.
$$\ell^*(f):=\sup\{a:f(a)=0\},\;\;\;\;\;\;\;\mbox{and}
\;\;\;\;\;\;\ell_*(f):=\inf\{a:f(a)=0\}.$$$\;$

\noindent If $f(s)f(s')<0$ for some $s,s'\in\R$, then clearly $\ell_*(f)$ and $\ell^*(f)$ are finite, and  $-\infty<\ell_*(f)<\ell^*(f)<\infty$. In fact, in this case, $f\equiv 0$ in the closed interval $[\ell_*(f),\ell^*(f)]$. Moreover, we note that $\ell_*(f)=\ell^*(f)$ if and only if $f(s)=0$ for exactly one $s\in\R$. If $\ell_*(f)=-\infty$, and $\ell^*(f)=\infty$, then $f\equiv0$ in $\R$.  Finally let us also note that if $\ell^*(f)=-\infty$ or $\ell_*(f)=\infty$, then $f$ never vanishes in $\R$.\\

\noindent In the following theorem we assume that at least one of $\ell_*(f)$ or $\ell^*(f)$ is finite, for otherwise the result is true.\\

\begin{thm}\label{cm} Let $f\in C(\R,\R)$ be non-decreasing. Let $u,v\in C(\overline{\Om})$ such that
$$\Df u\geq f(u)\;\mbox{in}\;\Om,\;\;\Df v\leq  f(v)\;\mbox{in}\;\Om.$$ Suppose $v\geq \ell^*(f)$ or $u\leq \ell_*(f)$ on $\p\Om$. If $u\leq v$ in $\p\Om$, then $u\leq v$ in $\Om$.

\end{thm}

\pf We first take up the case $v\geq \ell^*(f)$ on $\p\Om$. Then $\ell^*(f)<\infty$.  Unless noted otherwise, we use the notation $\ell^*:=\ell^*(f)$ and $\ell_*:=\ell_*(f)$ in the proof below. We start by observing that any solution $w\in C(\overline{\Om})$ of $\Df w\leq f(w)$ in $\Om$ with $w\geq \ell^*$ on $\p\Om$ satisfies $w\geq \ell^*$ in $\Om$.  To see this, suppose $w(x)<\ell^*$ for some $x\in\Om$. Consider the non-empty open set $\Om_0=\{x\in\Om: w(x)<\ell^*\}$. Then $\Df w\leq f(w)\leq f(\ell^*)=0$ in $\Om_0$, and $w=\ell^*$ on $\p\Om_0$. Thus, by the maximum principle, we see that $w\geq \ell^*$ in $\Om_0$, which is a contradiction. Thus the claim follows. Now suppose $u(x)>v(x)$ for some $x\in\Om$. Consider the non-empty open set $\Om_0=\{x\in\Om:u(x)>v(x)\}$. In $\Om_0$ we have
$$\Df u(x)\geq f(u(x)):=k(x).$$ Since $u(x)>\ell^*$ for all $x\in\Om_0$, it is clear that $k\in C(\overline{\Om}_0)$ is positive in $\Om_0$. On the other hand in $\Om_0$, we see that
$$\Df v(x)\leq f(v(x))\leq f(u(x))=k(x).$$ Since $u=v$ on $\p\Om_0$, we conclude that $u\leq v$ in $\Om_0$ by Lemma \ref{luw}. But again, this is a contradiction.\\

\noindent Now let us suppose that $u\leq \ell_*(f)$ on $\p\Om$. Then $\ell_*(f)>-\infty$. Let $h(t):=-f(-t)$, $\widehat{u}:=-u$ and $\widehat{v}:=-v$. Then we see that
$$\Df \widehat{v}\geq h(\widehat{v})\;\mbox{in}\;\Om,\;\;\Df \widehat{u}\leq h(\widehat{u})\;\mbox{in}\;\Om\;\;\;\mbox{and}\;\;\widehat{v}\leq \widehat{u}\;\;\mbox{on}\;\p\Om.$$
Furthermore, $\ell^*(h)=-\ell_*(f)$ and hence $\widehat{u}\geq \ell^*(h)$ on $\p\Om$. Hence by the above, we conclude that $\widehat{v}\leq \widehat{u}$ in $\Om$, that is $u\leq v$ in $\Om$.\epf

\begin{thm}\label{cm1} Let $f\in C(\R,\R)$ be non-decreasing. Let $u,v\in C(\overline{\Om})$ such that
\eqRef{fs}
\Df u\geq f(u)\;\;\mbox{in}\;\Om,\;\;\mbox{and}\;\;\;\Df v\leq f(v)\;\mbox{in}\;\Om.
\end{equation}
If $\sup_{\p\Om}u\leq \inf_{\p\Om}v$, then $u\leq v$ in $\Om$.

\end{thm}

\pf For ease of notation, set $\ell_*=\ell_*(f)$ and $\ell^*=\ell^*(f)$. If $\ell_*=\infty$ or $\ell^*=-\infty$, then $f>0$ in $\R$ or $f<0$ in $\R$. In this case, the result follows from \cite[Lemma 4.3]{BMO}. If $\ell_*=-\infty$ and $\ell^*=\infty$, then $f\equiv 0$ in $\R$. Since $u$ is infinity-subharmonic and $v$ is infinity super-harmonic, the result follows. So we assume that $\ell_*$ or $\ell^*$ is finite. That is we assume that $f$ is non-trivial and vanishes somewhere. Let us fix a constant $c$ such that $\sup_{\p\Om}u\leq c\leq\inf_{\p\Om}v$.
\begin{description}
\item[Case 1] $-\infty<\ell_*\leq\ell^*<\infty$. In the case when $c\le \ell_*$ or when $\ell^*\le c$, the desired result follows from Theorem \ref{cm}. Now suppose that $\ell_*<c<\ell^*$.  We claim that $\ell_*\leq v$ and $u\leq \ell^*$ in $\Om$. We will prove the first, the proof of the second being similar. Suppose $v(x)<\ell_*$ for some $x\in\Om$. Consider the set $\Om_0=\{x\in\Om:v(x)<\ell_*\}$. Then on $\Om_0$ we have $\Df v\leq f(v)\leq f(\ell_*)=0$, and $v=\ell_*$ on $\p\Om_0$. Then $v\geq \ell_*$ in $\Om_0$ by a standard comparison principle. But this is an obvious contradiction. Therefore the claim holds. Now if $u$ and $v$ satisfy (\ref{fs}) and $\Om_0:=\{x\in\Om:u(x)>v(x)\}$ is non-empty, then $\Df u\geq f(u)\geq f(v)\geq f(\ell_*)=0$ and $\Df v\leq f(v)\leq f(u)\leq f(\ell^*)=0$ in $\Om_0$. Since $u=v$ on $\p\Om_0$, we conclude that $u\leq v$ in $\Om_0$, which is a clear contradiction.\\

\item[Case 2] Suppose $\ell_*=-\infty$ and $\ell^*<\infty$. If $c\geq \ell^*$, then the desired inequality follows from Theorem \ref{cm}. So let us suppose that $c<\ell^*$. We note that $f(t)=0$ for all $t\leq \ell^*$. We show that $u\leq \ell^*$ in $\Om$. Suppose $\Om_0:=\{x\in\Om:u(x)>\ell^*\}$. Then $\Df u\geq f(u)\geq 0$ in $\Om_0$ and $u=\ell^*$ on $\p\Om_0$ implies that $u\leq \ell^*$ in $\Om_0$, and this is a contradiction. Now if $u$ and $v$ satisfy (\ref{fs}), and $\Om_0=\{x\in\Om:u(x)>v(x)\}$ is non-empty, then $\Df u\geq f(u)=0$ and $\Df v\leq f(v)\leq f(u)=0$ in $\Om_0$ and since $u=v$ on $\p\Om_0$, we conclude that $u\leq v$ in $\Om_0$. Again this is a contradiction.\\

\item[Case 3] Suppose $\ell_*>-\infty$ and $\ell^*=\infty$. Then we consider $h(t):=-f(-t)$, and $\widehat{u}:=-u$ and $\widehat{v}:=-v$. We note that $h(t)\geq 0$ for all $t\in\R$, and $h(-e)>0$. Furthermore, we see that $\sup_{\p \Om}\widehat{v}\leq \inf_{\p\Om} \widehat{u}$ and
    $$\Df \widehat{v}\geq h(\widehat{v})\;\;\mbox{in}\;\;\Om,\;\;\;\;\Df \widehat{u}\leq h(\widehat{u})\;\;\mbox{in}\;\;\Om.$$ Therefore the conclusion follows from Case 2.

\end{description}
\qed\\

\begin{cor}[Uniqueness] Let $f:\R\rar\R$ be non-decreasing. If $u\in C(\overline{\Om})$ and $v\in C(\overline{\Om})$ are solutions of (\ref{dp}) then
$$\|u-v\|_{\Om;\,\infty}\leq \mbox{Osc}_{\p\Om}b.$$
\end{cor}

\pf Let us note each of $u+\sup_{\p\Om} b-\inf_{\p\Om} b$ and $v+\sup_{\p\Om}b-\inf_{\p\Om}b$ is a super-solution of (\ref{dp}). Therefore, by the above theorem, we note that
see that $u+\inf_{\p\Om} b-\sup_{\p\Om}\leq v$ in $\Om$, and $v+\inf_{\p\Om}b-\sup_{\p\Om}b\leq u$ in $\Om$. The stated conclusion then follows from these two inequalities.\quad\\

\begin{cor}[Uniqueness] Let $f:\R\rar\R$ be non-decreasing, and $c$ a constant. Then the Dirichlet problem
$$\Df u=f(u)\;\;\mbox{in}\;\;\Om\;\;\;\;\mbox{and}\;\;\;u=c\;\;\mbox{on}\;\;\p\Om,$$ has a unique solution.
\end{cor}

\pf This follows from Theorem \ref{cm1}.
\qed\\

\section{Appendix}

\begin{center}
\NI{\bf Part I}
\end{center}

In Part I of this appendix we prove three lemmas that have been used in Section 4. See Theorem \ref{uss}, Lemma \ref{apb2}, and Remark \ref{eig2}.

Let $\ell\in\R$. We start with a continuous, non-decreasing function $h:[\ell,\infty)\rar[0,\infty)$ such that $h(t)>0$ for $t>\ell$. Let

$$\psi(t):=\int_t^a \frac{ds}{(H(a)-H(s))^{1/4}},\;\;\;\;\;\;
\ell\leq t\leq a\;\;\;\;\tx{where}\;\;\;\;H(s)=\int_\ell^s h(\tau)\,d\tau.$$$\;$\\
Hereafter, $\psi$ and $H$ will stand for the quantities defined above.\\

The first lemma provides upper and lower bounds for $\psi(t)$.

\begin{lem}\label{app2} Let $h:[\ell,\infty)\rar [0,\infty)$ be non-decreasing, and continuous such that $h(t)>0$ for $t>\ell$. Then
\eqRef{pss}\frac43\left(\frac{(a-t)^3}{h(a)}\right)^{1/4}\leq \psi(t)\leq \frac43\left(\frac{(a-t)^4}{H(a)}\right)^{1/4}\;\;\;\;\;\;\;\;
\;\;\;\ell\leq t\leq a.
\end{equation}
\end{lem}

\pf
To prove the first inequality, we note that
\ben
H(a)-H(s)=\int_a^s h(t)\;dt \leq (a-s)h(a).
\een
Therefore for $\ell\leq t\leq a$ we have
$$\int_t^a\frac{dt}{(H(a)-H(s))^{1/4}}\,ds\geq
\frac{1}{h(a)^{1/4}}\int_t^a(a-s)^{-1/4}\,dt=\frac43
\left(\frac{(a-t)^3}{h(a)}\right)^{1/4}.$$$\;$\\

\nin Next we proceed to prove the second inequality. Let us start by noting that for $\ell<t$, the function $H(t)$ is convex in $t$. This follows since $H'(t)=h(t)$ is non-decreasing. As $H(\ell)=0$, it follows that
\eqRef{ssss}
\frac{H(t)}{H(a)}\leq\frac {t-\ell}{a-\ell},\;\;\;\;\;\ell<t\leq a.
\end{equation}
Employing (\ref{ssss}), we see that
\begin{align}
\int_t^a(H(a)-H(\xi))^{-1/4}\,d\xi&=
H(a)^{-1/4}\int_t^a\left(1-\frac{H(\xi)}{H(a)}\right)^{-1/4}\,d\xi
\notag\\[.2cm]
&\leq
H(a)^{-1/4}\int_t^a\left(1-\frac {\xi-t}{a-t}\right)^{-1/4}\,d\xi=
\left(\frac{a-t}{H(a)}\right)^{1/4}\int_t^a(a-\xi)^{-1/4}\,d\xi
\notag\\[.2cm]
&=\frac43\left(\frac{(a-t)^4}{H(a)}\right)^{1/4}.\notag
\end{align}
\qed\\

Let us consider the following conditions.

\begin{align}
&\int_\ell^{\ell+1}\frac{1}{\sqrt[4]{H(t)}}\,dt<\infty,
\;\;\;\;\;\;\;\;\;\;\;\;\;\;\;\;\;\;\;\;\;\;\;\;\;\;\;\;\;\;\;\;\;\;
\;\;\;\;\;\;\;\;\;\;\;\;\;\label{ho}\\[.3cm]
&\lim_{t\rar\infty}\frac{(t-\ell)^3}{h(t)}=\infty\label{ht}
\end{align}$\;$\\

Recalling the definitions of $\zeta$ in (\ref{dd2}) (also see (\ref{zta})), we observe that
$$\zeta(a)=\psi(\ell).$$
We have the following result.\\

\begin{lem}\label{app3} Let $h:[\ell,\infty)\rar (0,\infty)$ be non-decreasing continuous, and $\zeta$ be as defined above. Furthermore, we suppose that $h(t)>0$ for $t>\ell$.
\begin{enumerate}
\item If $h$ satisfies condition (\ref{ho}), then $\lim_{a\rar \ell^{\small +}}\zeta(a)=0$.

\item If $h$ satisfies condition (\ref{ht}), then $\lim_{a\rar\infty}\zeta(a)=\infty$.
\end{enumerate}
\end{lem}

\pf The lemma follows from Inequality (\ref{pss}) on taking $t=\ell$ and combining with the observation that
$$\frac43\left(\frac{(a-\ell)^4}{H(a)}\right)^{1/4}\leq \int_\ell^a\frac{1}{\sqrt[4]{H(t)}}\,dt.$$
Consequently,
$$\frac43\left(\frac{(a-\ell)^3}{h(a)}\right)^{1/4}\leq \zeta(a)\leq \int_\ell^a\frac{1}{\sqrt[4]{H(t)}}\,dt.$$
Clearly, if (\ref{ho}) holds then conclusion (1) follows, and if (\ref{ht}) holds then it is easily seen that conclusion (2) holds.\epf

Finally we state the following lemma which was used in the proof of Lemma \ref{apb2}.

\begin{lem}\label{app4} Let $h:[\ell,\infty)\rar[0,\infty)$ such that $h(\ell)=0$ and $h(t)>0$ for $t>\ell$. Then given $a>0$, there exists an $R>0$ such that the following initial boundary value problem has a solution $\phi\in C^2((0,R])\cap C^1([0,R])$.

\eqRef{invp}
\left\{\begin{array}{ll}
         (\phi'(r))^2\phi''(r) =-h(\phi(r))\;\;&\tx{for}\;\;0<r<R \\[.2cm]
         \phi(0)=a,\;\; \phi'(0)=0,\;\;\;\phi(R)=\ell\;\;&
       \end{array}
\right.
\end{equation}
\end{lem}

\pf

\nin Given $a>\ell$ we define
$$\psi(t)=\frac{1}{\sqrt{2}}\int_t^a\frac{ds}{(H(a)-H(s))^{1/4}},\;\;\;\;\;
\ell<t<a.$$ Let $R:=\psi(\ell)$. Then $\psi:[\ell,a]\rar [0,R]$ is a continuous, decreasing function, and if $\phi$ is its inverse function then, $\phi\in C^2((0,R])\cap C^1([0,R])$ is also a decreasing function. Observing that
$$r=\frac{1}{\sqrt{2}}\int_{\phi(r)}^a\frac{ds}{(H(a)-H(s))^{1/4}},$$
direct computation shows that
\begin{align*}
\phi'(r)&=- [\;4\{H(a)-H(\phi(r))\}\;]^{1/4}&&\tx{and}\\[.2cm]
\phi''(r)&=-h(\phi(r))   [\;4\{H(a)-H(\phi(r))\}\;]^{-1/2}.&&
\end{align*}

\nin Therefore $\phi$ is a solution of (\ref{invp}) as claimed.\qed

\newpage

\begin{center}
\NI{\bf Part II.}
\end{center}

In this part, we present some remarks about the problem
$$\Delta_{\infty}u+\lm a(x)u^3=0\;\;\;\mbox{in $\Om$ and $u=0$ on $\partial \Om$},$$
where $\lm>0$ is a parameter and $a$ is a non-negative and non-trivial continuous function in $\Om$. The reader is directed to Case 2 of Remark \ref{examp} where it was shown that if
\eqRef{crtrn}
\lm <\frac{1}{\s^3 \sup_\Om a(x) R^4},
\end{equation}
then $u(x)\equiv 0$ is the only solution. Also see Remark \ref{eig2} in this context. Let $a\in C(\Om)\cap L^\infty(\Om)$ such that $a(x)\geq 0$ and $a(x)\not\equiv0$ in $\Om$. As discussed in Remark \ref{examp}, the existence of a non-trivial solution to the problem
\eqRef{egnv}
\Delta_{\infty}u+\lm a(x) u^3=0 \;\;\;\;\mbox{in $\Om$ and $u=0$ on $\partial \Om$}
\end{equation}
requires a description of the permissible set of $\lm$ (hereafter referred to as eigenvalues) and corresponding non-trivial solutions (referred to as eigenfunctions). The techniques developed in our current work do not readily apply to this question and, as of now, it is not clear how they may be adapted to this situation. A preliminary calculation of the Appendix of \cite{BMO} (see (8.5) therein) shows that, in the case of the ball, there is an eigenvalue $\lm>0$ for which $u>0$ exists. The analysis is done only for the radial case. It is well known that the set of eigenvalues is infinite and unbounded in the case of elliptic operators, such as the Laplacian or even the $p$-Laplacian, for the Dirichlet version of the problem. Since the $p$-Laplacian has close connections to the infinity-Laplacian, it should be stated that a complete characterization of the set of the eigenvalues for the $p$-Laplacian is not yet available, at least to our knowledge. However, it is well-known that an eigenfunction has one sign if and only if it is the first eigenfunction. As a matter of fact, this holds true for a large class of elliptic operators. It is in this connection that we state and prove the following result for (\ref{egnv}).   \\

\begin{lem}\label{egn1}
Let $\Om\su \R^N$ be a bounded domain, and $a(x)\in C(\Om)\cap L^{\infty}(\Om),\;a(x)\ge 0$ and $a(x)\not\equiv 0.$ Set $M=\sup_{\Om} a(x).$ Suppose that $(\lm,u)$ is an eigen-pair of (\ref{egnv}) then
$$\lm\ge \frac{1}{\s^3 M R^{4}},$$
where $R$ is the radius of the out-ball for $\Om$. For $0<\al<1$, define $\Om_{\al}=\{x\in\Om:\;a(x)>\al M\}$. In the event, $u>0$ then
$$0<\frac{1}{\s^3{M} R^{4}}\le \lm\le {4\left(\frac43\right)^3}\frac{1}{\s^3M}\inf_{0<\al\le 1}\left(\frac{1}{\al\rho_{\al}^4}\right)<\infty, $$
where $\rho_{\al}$ is the radius of the in-ball for $\Om_{\al}$.
Thus if an eigenvalue $\lm$ is large enough then the corresponding eigenfunction $u$ changes sign in $\Om$.
\end{lem}

\pf As discussed in Remark \ref{examp}, in order for Problem (\ref{egnv}) to have a non trivial solution $u\in C(\overline{\Om})$, the parameter $\lm$ should satisfy
$$\lm\ge \frac{1}{\s^3 M R^{4}},$$
where $R$ is the radius of the out-ball for $\Om$. \\

Let us assume that $(\lm,u)$, with $u\in C(\overline{\Om})$, satisfies
$$\Delta_{\infty}u+\lm a(x) u^3=0,\;\;u>0\;\;\;\mbox{in $\Om,\;\;\;$ and $u=0$ on $\partial \Om$}.$$
Since $\lm>0$ and $a(x)\ge 0$, it follows that $u$ is infinity super-harmonic in $\Om$ and satisfies the strong minimum principle. For $0<\al<1$, let $B_{\rho_{\al}}(z_{\al})$ be the in-ball for $\Om_{\al}$. For $0\le r\le \rho_{\al}$, define $m(r)=\inf _{\p B_r(z_{\al})}u$. Then $u\ge m(r)$ in $B_r(z_\al)$, $m(r)$ is decreasing and concave.

Consider the function
$$v(x)=m(z_{\al})\left(1-\frac{|x-z_{\al}|}{\rho_{\al}}\right),\;\;0<|x-z_\al|<\rho_{\al}.$$
It is clear that $v$ is infinity-harmonic in $B_{\rho_{\al}}(z_{\al})\setminus \{z_{\al}\}$. Since $u\ge 0$ on $\p B_{\rho_{\al}}(z_{\al})$ and $u(z_{\al})=m(z_{\al})$, by the Comparison Principle, Lemma \ref{luw}, $v\le u$ in $B_{\rho}(z_{\al})\setminus\{z_{\al}\}$. In particular, for $0\le \theta <1$,
\eqRef{mro}
\frac{m(z_{\al})}{m(\theta\rho_{\al})}\le \frac{1}{1-\theta}
\end{equation}
Next we consider the function
$$w(x)=\s\lm^{1/3}(\al M)^{1/3}m(\theta\rho_{\al})\left( \left(\theta\rho_{\al}\right)^{4/3}-|x-z_{\al}|^{4/3}\right)+m(\theta\rho_{\al}),\;\;\;\;0\le |x-z_{\al}|\le \theta\rho_{\al}.$$
A simple calculation shows that
$$\Delta_{\infty}w=-\al\lm M m(\theta\rho_{\al})^3\;\;\;\;\mbox{in $B_{\theta\rho_{\al}}(z_{\al})\;\;$ and $\;\;w=m(\theta\rho_{\al})$ on $|x-z_{\al}|=\theta\rho_{\al}$}.$$
Since $B_{\theta\rho_{\al}}(z_{\al})\subset \Om_{\al}$, $\Delta_{\infty}u=-\lm a(x) u^3<-\al\lm M m(\theta\rho_{\al})^3$, in $B_{\theta\rho_{\al}}(z_{\al})$, and $u\ge m(\theta\rho_{\al})$ on $|x-z_{\al}|=\theta\rho_{\al}$, the Comparison Principle, Lemma \ref{luw}, yields that
$w\le u$. Moreover,
$$w(z_{\al})=\s\lm^{1/3}(\al M)^{1/3}m(\theta\rho_{\al})\left(\theta\rho_{\al}\right)^{4/3}+m(\theta\rho_{\al})\le u(z_{\al})=m(z_{\al}).$$
Rearranging and using (\ref{mro}), we have
$$ \s\lm^{1/3}(\al M)^{1/3}\left(\theta\rho_\al\right)^{4/3}+1\le \frac{m(z_{\al})}{m(\theta\rho_{\al})}\le \frac{1}{1-\theta}.$$
Clearly,
$$\s \lm^{1/3}(\al M)^{1/3}\rho_\al^{4/3}\le \frac{1}{\theta^{1/3}(1-\theta)},\;\;\;0<\theta<1.$$
By computing the minimum of the right hand side, which occurs at $\theta=1/4$, we obtain
$$\lm\le{4\left(\frac43\right)^3\frac{1}{\s^3M\al \rho_{\al}^4}}.$$\epf

$\;$\\$\;$\\$\;$\\$\;$

\end{document}